\documentclass[]{elsart}
\usepackage{latexsym,amssymb,amsmath,amscd}

\parskip=\medskipamount

\newcommand{\A}{{\mathbb A}}
\newcommand{\F}{{\mathbb F}}
\newcommand{\PP}{{\mathbb P}}
\newcommand{\R}{{\mathbb R}}
\newcommand{\Z}{{\mathbb Z}}
\newcommand{\cA}{{\mathcal {A}}}
\newcommand{\cB}{{\mathcal {B}}}

\newcommand{\p}{{\partial}}

\begin{document}

\begin{frontmatter}

\title{On the calculation of {\rm UNil}$_*$}

\author{Frank Connolly}
\address{Department of Mathematics\\
University of Notre Dame\\
Notre Dame, IN 46556-5683\\
E-mail: connolly.1@nd.edu}
\author{Andrew Ranicki}
\address{School of Mathematics\\
 University of Edinburgh\\
 King's Buildings\\
 Edinburgh EH9 3JZ\\
 Scotland, UK\\
E-mail: a.ranicki@ed.ac.uk}

\begin{abstract}
Cappell's codimension 1 splitting obstruction surgery group
$\hbox{\rm UNil}_n$ is a direct summand of the Wall surgery
obstruction group of an amalgamated free product. For any ring
with involution $R$ we use the quadratic Poincar\'e cobordism
formulation of the $L$-groups to prove that
$$L_n(R[x])~=~L_n(R)\oplus \hbox{\rm UNil}_n(R;R,R)~.$$
We combine this with M. Weiss' universal chain bundle theory to
produce almost complete calculations of $\hbox{\rm
UNil}_*(\Z;\Z,\Z)$ and the Wall surgery obstruction groups
$L_*(\Z[D_{\infty}])$ of the infinite dihedral group
$D_{\infty}={\Z}_2*{\Z}_2$.
\end{abstract}
\end{frontmatter}

\maketitle


%%%%%%%%%%%%%%%%%%%%%%%%%%%%%%%%%%%%%%%%
\section*{Introduction}\label{sect1}

The nilpotent $K$- and $L$-groups of rings are a rich source of
algebraic invariants for geometric topology, giving results of two
types: if the groups are zero it is possible to solve the associated
splitting and classification problems, while if they are non-zero the
groups are infinitely generated and the solutions to the problems are
definitely obstructed.  See Bass \cite{B}, Farrell \cite{Fa1},
\cite{Fa2}, Farrell and Hsiang \cite{FH}, Cappell \cite{CaB},\cite{CaI},
Ranicki \cite{RES}, Connolly and Ko\'zniewski \cite{CK}.

The unitary nilpotent $L$-groups $\hbox{\rm UNil}_*$ arise as follows.
Suppose given a closed $n$-dimensional manifold $X$ which is expressed
as a union of codimension 0 submanifolds $X_1,X_{-1} \subseteq X$
$$X~=~X_1\cup X_{-1}$$
with
$$X_0~=~X_1\cap X_{-1}~=~\partial X_{-1}~=~\partial X_{1}\subseteq X$$
a codimension 1 submanifold. Assume $X,X_{-1},X_0,X_1$ are connected, and
that the maps $\pi_1(X_0)\to \pi_1(X_{\pm 1})$ are injective,
so that by the van Kampen theorem the fundamental group of $X$ is an
amalgamated free product
$$\pi_1(X)~=~\pi_1(X_1)*_{\pi_1(X_0)}\pi_1(X_{-1})$$
with $\pi_1(X_i) \to \pi_1(X)$ ($i=-1,0,1$) injective. Given another closed
$n$-dimensional manifold $M$ and a simple homotopy equivalence $f:M \to X$
there is a single obstruction
$$s(f) \in \hbox{\rm UNil}_{n+1}(R;\cB_1,\cB_{-1})$$
to deforming $f$ by an $h$-cobordism of
domains to a homotopy equivalence of the form
$$f_1\cup f_{-1}~:~M_1\cup M_{-1}~\to~ X_1\cup X_{-1}~$$
with $f_{\pm 1}:(M_{\pm 1},\partial M_{\pm 1})\to (X_{\pm
1},\partial X_{\pm 1})$ homotopy equivalences of manifolds with
boundary such that
$$f_1\vert~=~f_{-1}\vert~:~\partial M_1~ =~ \partial M_{-1} \to \partial X_1~=~\partial X_{-1}$$
and
$$R~=~\Z[\pi_1(X_0)]~,~
{\cB}_{\pm 1}~=~\Z[\pi_1(X_{\pm 1})\backslash \pi_1(X_0)]~.$$
Cappell \cite{CaB},\cite{CaI} proved geometrically that the free Wall
\cite{Wa} surgery obstruction groups $L_*=L^h_*$ of the fundamental
group ring
$$\Lambda~=~\Z[\pi_1(X)]~=~\Z[\pi_1(X_1)]*_{\Z[\pi_1(X_0)]}\Z[\pi_1(X_{-1})]$$
have direct sum decompositions
$$L_*(\Lambda)~=~\hbox{\rm UNil}_*(R;\cB_1,\cB_{-1}) \oplus
L'_*(\Z [\pi_1(X_0)] \to \Z [\pi_1(X_1)]\times \Z [\pi_1(X_{-1})])$$
with $L'_*$ appropriately decorated intermediate relative $L$-groups.
The split monomorphism
$$\hbox{\rm UNil}_{n+1}(R;\cB_1,\cB_{-1})~\to~L_{n+1}(\Lambda)~;~s(f)
\mapsto \sigma(g)$$ sends the splitting obstruction $s(f)$ to the
surgery obstruction $\sigma(g)$ of the `unitary nilpotent
cobordism' of \cite{CaI}, an $(n+1)$-dimensional normal map
cobordism between $f$ and a split homotopy equivalence. The
4-periodicity $L_*(\Lambda)=L_{*+4}(\Lambda)$ extends to a
4-periodicity
$$\hbox{\rm UNil}_*(R;\cB_1,\cB_{-1})~=~\hbox{\rm
UNil}_{*+4}(R;\cB_1,\cB_{-1})~.$$
\indent Farrell \cite{Fa2} obtained a remarkable factorization
$$\hbox{\rm UNil}_{n+1}(R;\cB_1,\cB_{-1})~\to~
\hbox{\rm UNil}_{n+1}(\Lambda;\Lambda,\Lambda)~\to~ L_{n+1}(\Lambda)$$
For this reason (and some others too) the groups $\hbox{\rm UNil}_*(R;R,R)$
for any ring with involution $R$ are of especial significance to us,
and we introduce the abbreviation:
$$\hbox{\rm UNil}_n(R)~ =~\hbox{\rm UNil}_n(R;R,R)~.$$
\indent But even the groups $\hbox{\rm UNil}_*(\Z)$ have remained
opaque for the last 30 years.  Cappell \cite{CaS},\cite{CaB},\cite{CaM}
proved that $\hbox{\rm UNil}_{4k}(\Z)=0$ and that $\hbox{\rm
UNil}_{4k+2}(\Z)$ is infinitely generated.  The $\hbox{\rm
UNil}$-groups $\hbox{\rm UNil}_*(R;\cB_1,\cB_{-1})$ are 2-primary
torsion groups.  Farrell \cite{Fa2} proved that $4\hbox{\rm UNil}_*(R)
= 0$, for any ring $R$.  Connolly and Ko\'zniewski \cite{CK} obtained
an isomorphism
$$\hbox{\rm UNil}_{4k+2}(\Z)~\cong~\bigoplus^{\infty}_1\F_2~,$$
together with
information on $\hbox{\rm UNil}_{4k+2}(R)$ for various Dedekind domains
and division rings.  But that is nearly all that is known.

The infinite dihedral group is a free product of two copies of the cyclic
group ${\Z}_2$ of order 2
$$D_{\infty}~=~{\Z}_2*{\Z}_2~.$$
Since the surgery obstruction groups $L_*(R[D_{\infty}])$ are hard to
compute directly, the split monomorphisms $\hbox{\rm UNil}_*(R)\to
L_*(R[D_{\infty}])$ are more useful in computing $L_*(R[D_{\infty}])$
from $\hbox{\rm UNil}_*(R)$ than the other way round.  Connolly and
Ko\'zniewski \cite{CK} expressed $\hbox{\rm UNil}_*(R;\cB_1,\cB_{-1})$
as the $L$-groups $L_*(\A_\alpha [x])$ of an {\it additive category}
with involution $\A_\alpha [x]$.  Although this expression did give new
computations of $\hbox{\rm UNil}_*(R)$, the $L$-theory of additive
categories with involution (Ranicki \cite{RAD}) is not in general very
computable.

The first goal of this paper therefore, is to provide a new description
for $\hbox{\rm UNil}_n(R)$ in terms of $L$-groups, which can be used to
computational advantage.
Cappell and Farrell observed that the infinite dihedral group
$D_{\infty}=\Z_2*\Z_2$ can also be viewed as an extension of $\Z$ by $\Z_2$
$$\{1\} \to \Z \to D_{\infty} \to \Z_2 \to \{1\}~,$$
so that the classifying space can be viewed both as a one-point union
$$K(D_{\infty},1)~=~K(\Z_2,1)\vee K(\Z_2,1)~=~
\R\PP^{\infty} \vee \R \PP^{\infty}$$
and as the total space of a fibration
$$K(\Z,1)~=~S^1 \to K(D_{\infty},1) \to K(\Z_2,1)~=~\R \PP^{\infty}~,$$
and that this should have implications for codimension 1 surgery
obstruction theory with $\pi_1=D_{\infty}$.  This observation was used
in Ranicki \cite{RES} (pp.  737--745) to prove geometrically that
for the group ring $R=\Z[\pi]$ of a finitely presented group $\pi$
$$\hbox{\rm UNil}_*(R)~=~NL_*(R)~=~\hbox{\rm ker}(L_*(R[x]) \to L_*(R))$$
with the involution on $R$ extended to $R[x]$ by $\overline{x}=x$,
and $R[x] \to R;x \mapsto 0$ the augmentation map.
The $NL$-groups are $L$-theoretic analogues of the nilpotent $K$-group
$$NK_1(R)~=~\hbox{\rm ker}(K_1(R[x]) \to K_1(R))~=~\widetilde{\hbox{\rm
Nil}}_0(R)$$
of Chapter XII of Bass \cite{B}, which is such that
$$K_1(R[x])~=~K_1(R)\oplus NK_1(R)~.$$

\noindent{\bf Theorem A.}
{\it For any ring with involution $R$
$$\hbox{\rm UNil}_*(R)~=~NL_*(R)$$
so that
$$L_n(R[x])~=~L_n(R) \oplus \hbox{\rm UNil}_n(R)~.$$}

We develop a new method for calculating $\hbox{\rm UNil}_*(R)$,
adopting the following strategy.  The symmetric $L$-groups
$L^*(R)$ of a ring $R$ with an involution $R \to R;x \mapsto
\overline{x}$ were defined by Mishchenko \cite{M} and Ranicki
\cite{RATI,RATII} to be the cobordism groups of symmetric
Poincar\'e complexes over $R$.  The quadratic $L$-groups $L_*(R)$
were expressed in \cite{RATI,RATII} as the cobordism groups of
quadratic Poincar\'e complexes over $R$, and the two types of
$L$-groups were related by an exact sequence
$$\dots \to L_n(R) \to L^n(R) \to \widehat{L}^n(R) \to L_{n-1}(R) \to \dots$$
with the hyperquadratic $L$-groups $\widehat{L}^*(R)$ the
cobordism groups of (symmetric,quadratic) Poincar\'e pairs. The
symmetric and hyperquadratic $L$-groups are not 4-periodic in
general, but there are defined natural maps
$$L^n(R) \to L^{n+4}(R)~,~\widehat{L}^n(R) \to \widehat{L}^{n+4}(R)$$
(which are isomorphisms for certain $R$, e.g. a Dedekind ring or the
polynomial extension of a Dedekind ring). The 4-periodic versions of the
symmetric and hyperquadratic $L$-groups
$$L^{n+4*}(R)~=~\lim\limits_{k\to\infty}L^{n+4k}(R)~,~
\widehat{L}^{n+4*}(R)~=~\lim\limits_{k\to\infty}\widehat{L}^{n+4k}(R)$$
are related by an exact sequence
$$\dots \to L_n(R) \to L^{n+4*}(R) \to \widehat{L}^{n+4*}(R)
\to L_{n-1}(R) \to \dots~.$$ The theory of Weiss \cite{WeI,WeII}
identified $\widehat{L}^{n+4*}(R)$ with the `twisted $Q$-group'
$Q_n(B^R,\beta^R)$ of the `universal chain bundle' $(B^R,\beta^R)$
over $R$, which can be computed (more or less effectively) from
the Tate $\Z_2$-cohomology groups of the involution on $R$
$$\begin{array}{ll}
H_n(B^R)&=~\widehat{H}^n(\Z_2;R)\\[1ex]
&=~\{a \in R\,\vert\,\overline{a}=(-1)^na\}/\{b+(-1)^n\overline{b}\,\vert\,
b \in R\}~.
\end{array}$$
In Proposition \ref{2.1} we show that for a Dedekind ring with involution $R$
$$L^n(R[x])~=~L^n(R)~~,~~NL^n(R)~=~0$$
making the $\hbox{\rm UNil}$-groups
$$\hbox{\rm UNil}_n(R)~=~\hbox{\rm ker}(Q_n(B^{R[x]},\beta^{R[x]})\to Q_n(B^R,\beta^R))$$
accessible to computation.

\noindent{\bf Theorem B.} {\it For the ring $\Z$, we have:
$$\hbox{\rm UNil}_0(\Z)~=~0~,~\hbox{\rm UNil} _1(\Z)~=~0$$
and there is an exact sequence:
$$0\to \F_2[x]/\F_2 \overset{\psi^2-1}{\longrightarrow}
\F_2[x]/\F_2 \to \hbox{\rm UNil}_2(\Z)\to 0$$
with
$$\psi^2~:~\F_2[x] \to \F_2[x]~;~a \mapsto a^2$$
the Frobenius map.  $\hbox{\rm UNil}_3(\Z)$ is not finitely generated,
with $4 \hbox{\rm UNil}_3(\Z) = 0$.}

We now give an outline of the rest of this paper.

In \S1 we define the groups $\hbox{\rm UNil}_n(R)$ and the map
$c:\hbox{\rm UNil}_n(R)\to L_n(R[x])$, as well as the various
other morphisms and groups with which we will be working. Theorem
A is proved in \S1.

In \S2 we relate $\hbox{\rm UNil}_n(R)$ for Dedekind $R$ to the
group of symmetric structures on the universal chain bundle of
Weiss.  We then make the calculations necessary to prove Theorem
B.

We are grateful to the referee for various helpful suggestions.

This research was started during the visit of FXC to Edinburgh in
1994, which was supported by the SERC Visiting Fellowship
GR/J94891.

\section{ Fundamental Concepts.  The proof of Theorem
A.}\label{sect2}

\subsection{ Algebraic $L$-groups}

Throughout this paper $R$ denotes a ring with an involution
$$R \to R~;~r \mapsto \overline{r}~.$$
An $R$-module is understood to be a left $R$-module, unless a
right $R$-module action is specified. Given an $R$-module $P$ let
$P^t$ be the right $R$-module with the same additive group and
$$P^t \times R \to P^t~;~(x,r) \mapsto \overline{r}x~.$$
The dual of an $R$-module $P$ is the $R$-module
$$\begin{array}{l}
P^*~=~\hbox{\rm \hbox{\rm Hom}}_R(P,R)~,\\[1ex]
R \times P^* \to P^*~;~(r,f) \mapsto (x \mapsto f(x) \bar{r})~.
\end{array}$$
Write the evaluation pairing as
$$\langle~,~\rangle~:~ P^* \times P  \to R~;~
(f, x) \mapsto \langle f,x \rangle~=~f(x)~.$$
An element $\phi \in \hbox{\rm Hom}_R(P,P^*) $ determines a
sesquilinear form on $P $
$$\langle~,~\rangle_\phi~:~P \times P \to R~;~(x,y) \mapsto \langle
\phi(x),y \rangle~,$$
and we identify $\hbox{\rm Hom}_R(P,P^*) $ with the additive group of
such forms.  The dual of a f.g.  (= finitely generated) projective
$R$-module $P$ is a f.g.  projective $R$-module $P^*$, and the morphism
$$P \to P^{**}~;~x \mapsto (f \mapsto \overline{f(x)})$$
is an isomorphism, which we shall use to identify
$$P^{**}~=~P~,$$
and to define the $\epsilon$-duality involution
$$T_\epsilon~:~\hbox{\rm Hom}_R(P,P^*) \to \hbox{\rm Hom}_R(P,P^*)~;~
\phi \mapsto \epsilon \phi^*~~,~~\langle x,y
\rangle_{\phi^*}~=~\overline{\langle y,x \rangle}_\phi~.$$

For $\epsilon=\pm 1$, any $R$-module chain complex $C$ and any
$\Z[\Z_2]$-module chain complex $X$ define the $\Z$-module chain
complexes
$$\begin{array}{l}
X^{\%}(C,\epsilon)~=~\hbox{\rm Hom}_{\Z[\Z_2]}(X,C^t\otimes_RC)~,\\[1ex]
X_{\%}(C,\epsilon)~=~X\otimes_{\Z[\Z_2]}(C^t\otimes_RC)
\end{array}$$
with $T \in \Z_2$ acting on $C^t\otimes_RC$ by the signed
transposition isomorphisms
$$T_{\epsilon}~:~C_p^t\otimes_R C_q \to C_q^t\otimes_R C_p~;~
    x \otimes y \mapsto (-1)^{pq}\epsilon y \otimes x~.$$
We shall be mainly concerned with finite chain complexes $C$ of
f.g. projective $R$-modules, in which case we identify
$$C^t\otimes_RC~=~{\rm Hom}_R(C^*,C)$$
using the natural $\Z$-module isomorphisms
$$C_p^t \otimes_R C_q \to {\rm Hom}_R(C^p,C_q)~;~
    x \otimes y \mapsto (f \mapsto \overline{f(x)}.y)$$
with $C^p=(C_p)^*$. The signed transposition isomorphisms
correspond to the signed duality isomorphisms
$$T_{\epsilon}~:~{\rm Hom}_R(C^p,C_q) \to {\rm Hom}_R(C^q,C_p)~;~
    \phi \mapsto (-1)^{pq} \epsilon \phi^*~.$$
\indent As in Ranicki \cite{RATI,RATII} the group of
$n$-dimensional $\epsilon$-symmetric (resp.
$\epsilon$-hyperquadratic, resp. $\epsilon$-quadratic) structures
on $C$ is defined by:
$$\begin{array}{l}
Q^n(C,\epsilon)~=~H_n(W^\%C)~~,~~\widehat Q^n(C,\epsilon)~=~
H_n(\widehat W^\%C)~,\\[1ex]
Q_n(C,\epsilon)~=~H_n(W_{\%}C)~=~H_n((W^{-*})^\%C)
\end{array}$$
where $W$ (resp. $ \widehat W$) denotes the standard
free $ \Z[\Z_2]$-module resolution of $\Z$ (resp. complete
resolution) and
$$W^{-*}~=~\hbox{\rm Hom}_{\Z[\Z_2]}(W,\Z[\Z_2])~.$$
If $S^{-1}W^{-*}$ denotes the desuspension of $W^{-*}$, the short exact
sequence
$$ 0\to S^{-1}W^{-*}\to \widehat W\to W\to 0 $$
induces the  exact sequence:
\begin{equation}\label{3a}
\dots \to Q_n(C,\epsilon)\to Q^n(C,\epsilon)\overset{J}{\to } \widehat
Q^n(C,\epsilon)\to Q_{n-1}(C,\epsilon)\to \dots~.
\end{equation}

Given a f.g. projective $R$-module $P$ define the 0-dimensional
f.g. projective $R$-module chain complex
    $$C~:~\dots \to 0 \to C_0~=~P^* \to 0 \to \dots~.$$
At the risk of notational confusion, the 0-dimensional
$\epsilon$-symmetric and $\epsilon$-quadratic $Q$-groups of $C$
are written
$$\begin{array}{l}
Q^0(C,\epsilon)~=~Q^{\epsilon}(P)~=~\hbox{\rm
ker}(1-T_\epsilon:\hbox{\rm Hom}_R(P,P^*) \to \hbox{\rm
Hom}_R(P,P^*))~,\\[1ex]
Q_0(C,\epsilon)~=~Q_{\epsilon}(P)~=~\hbox{\rm
coker}(1-T_\epsilon:\hbox{\rm Hom}_R(P,P^*) \to \hbox{\rm
Hom}_R(P,P^*))~.
\end{array}$$

\begin{defn} {\rm
An {\it $\epsilon$-symmetric form} $(P,\phi)$ (resp. an {\it
$\epsilon$-quadratic form} $(P,\psi)$) over $R$ is a f.g.
projective $R$-module $P$ together with an element $\phi \in
Q^{\epsilon}(P)$ (resp. $\psi \in Q_{\epsilon}(P)$). The form is
{\it nonsingular} if the $R$-module morphism
$$\phi:P \to P^*~\hbox{(resp.  }
N_{\epsilon}(\psi)=(1+T_\epsilon)\psi:P \to P^*)$$
is an isomorphism.}\hfill$\qed$
\end{defn}
\indent We refer to Ranicki
\cite{RATI,RATII},\cite{RES},\cite{RAP} for various accounts of
the construction of the free $\epsilon$-symmetric (resp.
quadratic) $L$-groups $L^n(R,\epsilon)$ (resp. $L_n(R,\epsilon)$)
as the cobordism groups of $n$-dimensional $\epsilon$-symmetric
(resp. $\epsilon$-quadratic) Poincar\'e complexes over $R$
$(C,\phi \in Q^n(C,\epsilon))$ (resp. $(C,\psi \in
Q_n(C,\epsilon)))$ with
$$C~:~\dots \to 0 \to C_n \to C_{n-1} \to \dots \to C_0\to 0 \to \dots$$
an $n$-dimensional f.g. free $R$-module chain complex.
The projective $L$-groups $L^*_p(R,\epsilon)$ (resp. $L^p_*(R,\epsilon)$)
are constructed in the same way, using f.g. projective $C$.

The suspension of an $R$-module chain complex $C$ is the $R$-module
chain complex $SC$ with
$$d_{SC}~=~d_C~:~(SC)_{r+1}~=~C_r \to (SC)_r~=~C_{r-1}~.$$
As in Ranicki \cite{RATI,RATII} (p. 105) use the natural
$\Z$-module isomorphisms
$$S^2(W^{\%}(C,\epsilon))~\cong~W^{\%}(SC,-\epsilon)~~,~~
S^2(W_{\%}(C,\epsilon))~\cong~W_{\%}(SC,-\epsilon)$$
to identify
$$Q^n(C,\epsilon)~=~Q^{n+2}(SC,-\epsilon)~~,~~
Q_n(C,\epsilon)~=~Q_{n+2}(SC,-\epsilon)$$
and to define the skew-suspension maps
$$\begin{array}{l}
\overline{S}^n~:~L^n(R,\epsilon) \to L^{n+2}(R,-\epsilon)~;~
(C,\phi) \mapsto (SC,\phi)~,\\[1ex]
\overline{S}_n~:~L_n(R,\epsilon) \to L_{n+2}(R,-\epsilon)~;~(C,\psi)
\mapsto (SC,\psi)~.
\end{array}$$

\begin{defn} \label{1-dim}
{\rm A ring $R$ is {\it $1$-dimensional} if it is hereditary
and noetherian, or equivalently if every submodule of a f.g.
projective $R$-module is f.g.  projective.}
\vskip-5mm
$$\eqno{\qed}$$
\end{defn}
\vskip-4mm

In particular, Dedekind rings are $1$-dimensional.

\begin{prop} \label{below} {\rm (\cite{RATI,RATII})}\\
{\rm (i)} For every ring with involution $R$ the
$\pm\epsilon$-quadratic skew-suspension maps $\overline{S}_n$ are
isomorphisms, so that
$$L_n(R,\epsilon)~=~L_{n+2}(R,-\epsilon)~=~L_{n+4}(R,\epsilon)~,$$
with $L_{2n}(R,\epsilon)=L_0(R,(-1)^n\epsilon)$ the Witt group of
stable isometry classes of nonsingular $(-1)^n\epsilon$-quadratic forms
over $R$.  \\
{\rm (ii)} If $R$ is $1$-dimensional then the $\pm\epsilon$-symmetric
skew-suspension maps $\overline{S}^n$ are isomorphisms, so that
$$L^n(R,\epsilon)~=~L^{n+2}(R,-\epsilon)~=~L^{n+4}(R,\epsilon)~,$$
with $L^{2n}(R,\epsilon)=L^0(R,(-1)^n\epsilon)$ the Witt group of
stable isometry classes of nonsingular $(-1)^n\epsilon$-symmetric forms
over $R$.
\end{prop}
\noindent{\it Proof.}{\rm By algebraic surgery below the middle
dimension, given by Proposition I.4.3 of \cite{RATI,RATII} for
(i), and Proposition I.4.5 of \cite{RATI,RATII} for (ii).}
\hfill\qed

For $\epsilon=1$ we write
$$\begin{array}{l}
X^{\%}(C,1)~=~X^{\%}C~~,~~X_{\%}(C,1)~=~X_{\%}C~~,\\[1ex]
Q^*(C,1)~=~Q^*(C)~~,~~\widehat{Q}^*(C,1)~=~\widehat{Q}^*(C)~~,~~
Q_*(C,1)~=~Q_*(C)~~,\\[1ex]
L^*(R,1)~=~L^*(R)~~,~~L_*(R,1)~=~L_*(R)~~.
\end{array}$$
The hyperquadratic $Q$-groups $\widehat{Q}^*(C)$ are used in
Section \ref{sect3} to define chain bundles.

\subsection{The nilpotent $L$-groups $L\hbox{\rm Nil}$,
$L\widetilde{\hbox{\rm Nil}}$}

Theorem A identifies the unitary nilpotent $L$-groups $\hbox{\rm
UNil}_*(R)$ with the nilpotent $L$-groups $L\widetilde{\hbox{\rm
Nil}}_*(R)$, whose definition we now recall.

We start with nilpotent $K$-theory.

\begin{defn} {\rm (i) An {\it $R$-nilmodule} $(P,\nu)$ is a
f.g. projective $R$-module $P$ together with a nilpotent
endomorphism $\nu:P \to P$, so that
$$\nu^N~=~0~:~P \to P$$
for some $N \geqslant 1$.\\
(ii) A {\it morphism} of $R$-nilmodules $f:(P,\nu) \to (P',\nu')$ is an
$R$-module morphism $f:P \to P'$ such that $\nu'f=f\nu :P \to P'$.\\
(iii) The {\it nilpotent $K$-groups} of $R$ are defined to be the $K$-groups
$$\hbox{\rm Nil}_*(R)~=~K_*(\hbox{\rm Nil}(R))$$
of the exact category $\hbox{\rm Nil}(R)$ be of $R$-nilmodules.
The {\it reduced nilpotent $K$-groups}
$$\widetilde{\hbox{\rm Nil}}_*(R)~=~\hbox{\rm ker}(\hbox{\rm Nil}_*(R)
\to K_*(R))$$
are such that
$$\hbox{\rm Nil}_*(R)~=~K_*(R)\oplus \widetilde{\hbox{\rm Nil}}_*(R)~.$$
{\rm (iv)} The {\it $NK$-groups} of $R$ are defined by
$$NK_*(R)~=~\hbox{\rm ker}(K_*(R[x]) \to K_*(R))~,$$
so that
$$K_*(R[x])~=~K_*(R) \oplus NK_*(R)~.\eqno{\qed}$$}
\end{defn}

\begin{prop} {\rm (Bass \cite{B})}\\
{\rm (i)} There is a natural identification
$$NK_1(R)~=~\widetilde{\hbox{\rm Nil}}_0(R)$$
using the split injection
$$\widetilde{\hbox{\rm Nil}}_0(R) \to K_1(R[x])~;~
(P,\nu) \mapsto \tau(1+x\nu:P[x] \to P[x])~.$$
{\rm (ii)} If $R$ is $1$-dimensional then
$$\widetilde{\hbox{\rm Nil}}_0(R)~=~0~.$$
\end{prop}
\noindent{\it Proof.} (i) See Chapter XII of \cite{B}.\\
(ii) Given a nilmodule $(P,\nu)$ with $\nu^N=0:P \to P$ for some
$N \geqslant 1$ define the nilmodules
$$(P',\nu') ~=~(\hbox{\rm ker}(\nu),0)~~,~~
(P'',\nu'')~=~(\hbox{\rm im}(\nu),\nu\vert)~,$$
using the $1$-dimensionality of $R$ to ensure that the $R$-modules
$\hbox{\rm ker}(\nu),\hbox{\rm im}(\nu) \subseteq P$ are f.g. projective.
It follows from the exact sequence
$$0 \to (P',\nu')  \to (P,\nu) \to (P'',\nu'') \to 0$$
that
$$[P,\nu]~=~[P',\nu']+[P'',\nu''] \in \hbox{\rm \hbox{\rm Nil}}_0(R)~.$$
Now $\nu'=0$, $(\nu'')^{N-1}=0$, so proceeding inductively we obtain
$$[P,\nu]~=~\sum\limits^N_{i=1}[\hbox{\rm ker}(\nu^i)/\hbox{\rm
ker}(\nu^{i-1}),0]
\in K_0(R) \subseteq \hbox{\rm \hbox{\rm Nil}}_0(R)$$
and hence that $\widetilde{\hbox{\rm \hbox{\rm Nil}}}_0(R)=0$.
\hfill\qed

\begin{defn} {\rm An $n$-dimensional {\it $R$-nilcomplex}
$(C,\nu)$ is a $n$-dimensional f.g.  projective $R$-module chain complex
$$C~:~\dots \to 0 \to C_n \to C_{n-1} \to \dots \to C_1 \to C_0$$
together with a chain map $\nu:C \to C$ which is chain homotopy
nilpotent, i.e.  such that $\nu^N \simeq 0:C \to C$ for some integer
$N\geqslant 1$.\hfil\qed}
\end{defn}

\begin{prop} \label{2.13}
The chain equivalence classes of the following types of
chain complexes are in one-one correspondence:
\begin{itemize}
\item[\rm (i)] $n$-dimensional chain complexes of $R$-nilmodules
$$(C,\nu)~:~\dots \to 0 \to (C_n,\nu) \to (C_{n-1},\nu) \to \dots \to
(C_1,\nu) \to (C_0,\nu)~,$$
\item[\rm (ii)] $n$-dimensional $R$-nilcomplexes $(C,\nu)$,
\item[\rm (iii)] $(n+1)$-dimensional f.g.  projective $R[x]$-module
chain complexes
$$D~:~\dots \to 0 \to D_{n+1} \to D_n \to \dots \to D_1 \to D_0$$
such that
$$H_*(R[x,x^{-1}]\otimes_{R[x]}D)~=~0~.$$
\end{itemize}
\end{prop}
\noindent{\it Proof.} (i) $\Longrightarrow$ (ii)
An $n$-dimensional chain complex of $R$-nilmodules
is an $n$-dimensional $R$-nilcomplex.\\
(ii) $\Longrightarrow$ (iii) Given an $n$-dimensional $R$-nilcomplex
$(C,\nu)$ define the $(n+1)$-dimensional f.g.  projective $R[x]$-module
chain complexes
$$D~=~{\mathcal C}(x-\nu:C[x] \to C[x])$$
such that
$$\begin{array}{l}
H_*(R[x,x^{-1}]\otimes_{R[x]}D)~=~0~~,\\[1ex]
x~=~\nu~:~H_*(D)~=~H_*(C) \to H_*(D)~=~H_*(C)~.
\end{array}$$
(i) $\Longleftrightarrow$ (iii)
See Proposition 3.1.2 of Ranicki \cite{RES}.
\hfill\qed

In particular, it follows from Proposition \ref{2.13} that every
$n$-dimensional $R$-nilcomplex is chain equivalent to an
$n$-dimensional $R$-nilcomplex $(C,\nu)$ with $\nu^N=0:C \to C$ for
some $N \geqslant 1$ (rather than just $\nu^N \simeq 0$).

Now for nilpotent $L$-theory.

\begin{defn}\label{2.14}
{\rm (Ranicki \cite{RES}, p.\ 440, \cite{RHK} p.\ 470)\\
(i) The {\it $\epsilon$-symmetric $Q\hbox{\rm Nil}$-groups}
$Q\hbox{\rm Nil}^*(C,\nu,\epsilon)$ of an $R$-nilcomplex
$(C,\nu)$ are the relative $Q$-groups in the exact sequence
$$\dots \to Q^{n+1}(C,-\epsilon) \to Q\hbox{\rm Nil}^n(C,\nu,\epsilon) \to
Q^n(C,\epsilon) \overset{\Gamma_{\nu}}{\to} Q^n(C,-\epsilon) \to \dots$$
with
$$\Gamma_\nu~:~W^{\%}(C,\epsilon) \to W^{\%}(C,-\epsilon)~;~
\phi \mapsto (1\otimes \nu)\phi - \phi(\nu \otimes 1)~.$$
Similarly for the {\it $\epsilon$-quadratic $Q\hbox{\rm Nil}$-groups}
$Q\hbox{\rm Nil}_*(C,\nu,\epsilon)$, with an exact sequence
$$\dots \to Q_{n+1}(C,-\epsilon) \to
Q\hbox{\rm Nil}_n(C,\nu,\epsilon) \to Q_n(C,\epsilon)
\overset{\Gamma_{\nu}}{\to} Q_n(C,-\epsilon)\to \dots~.$$
(ii) An {\it $n$-dimensional $\epsilon$-symmetric Poincar\'e
nilcomplex over $R$} $(C,\nu,\delta\phi,\phi)$ is an
$n$-dimensional $R$-nilcomplex $(C,\nu)$ together with an element
$$(\delta\phi,\phi) \in Q\hbox{\rm Nil}^n(C,\nu,\epsilon)$$
such that $(C,\phi \in Q^n(C,\epsilon))$ is an
$n$-dimensional $\epsilon$-symmetric Poincar\'e complex over $R$.
The {\it $\epsilon$-symmetric $L\hbox{\rm Nil}$-group} $L\hbox{\rm
Nil}^n(R,\epsilon)$ is the cobordism group of $n$-dimensional
$\epsilon$-symmetric Poincar\'e nilcomplexes over $R$.
Similarly in the $\epsilon$-quadratic case, with $L\hbox{\rm
Nil}_n(R,\epsilon)$.\\
(iii) The {\it reduced $\epsilon$-symmetric $L\hbox{\rm Nil}$-groups}
are defined by
$$L\widetilde{\hbox{\rm Nil}}^*(R,\epsilon)~=~\hbox{\rm ker}(L\hbox{\rm
Nil}^*(R,\epsilon)\to L^*_p(R,\epsilon))~,$$
with
$$L\hbox{\rm Nil}^*(R,\epsilon)~=~L^*_p(R,\epsilon)\oplus
L\widetilde{\hbox{\rm Nil}}^*(R,\epsilon)~.$$
Similarly in the $\epsilon$-quadratic case, with $L\widetilde{\hbox{\rm
Nil}}_*(R,\epsilon)$.\\
(iv) Extend the involution to $R[x]$ by $\overline{x}=x$. Use the
augmentation map
$$R[x] \to R~;~x \mapsto 0$$
to define the {\it nilpotent $\epsilon$-symmetric $L$-groups of $R$}
$$NL^*(R,\epsilon)~=~\hbox{\rm ker}(L^*(R[x],\epsilon) \to L^*(R,\epsilon))$$
with
$$L^*(R[x],\epsilon)~=~L^*(R,\epsilon) \oplus NL^*(R,\epsilon)~.$$
Similarly for the {\it nilpotent $\epsilon$-quadratic $L$-groups}
$NL_*(R,\epsilon)$.\hfill\qed}
\end{defn}

\begin{prop}\label{2.20}\label{3g}
{\rm (i)} The $Q\hbox{\rm Nil}$-groups of an $R$-nilcomplex $(C,\nu)$
are the $Q$-groups of the $R[x,x^{-1}]$-contractible
f.g. projective $R[x]$-module chain complex
$$D~=~{\mathcal C}(x-\nu:C[x] \to C[x])$$
with
$$\begin{array}{l}
x~=~\nu~:~H_*(D)~=~H_*(C) \to H_*(D)~=~H_*(C)~,\\[1ex]
Q\hbox{\rm Nil}^n(C,\nu,\epsilon)~=~Q^{n+1}(D,-\epsilon)~,\\[1ex]
Q\hbox{\rm Nil}_n(C,\nu,\epsilon)~=~Q_{n+1}(D,-\epsilon)~.
\end{array}$$
An element $(\delta\phi,\phi) \in Q\hbox{\rm Nil}^n(C,\nu,\epsilon)$
corresponds to an element
$$\Phi \in Q^{n+1}(D,-\epsilon)~=~Q^{n-1}(S^{-1}D,\epsilon)$$
with
$$\phi_0~=~\Phi_0~:~H^{n+1-*}(D)~=~H^{n-*}(C) \to H_*(D)~=~H_*(C)~,$$
so that $(C,\nu,\delta\phi,\phi)$ is an $\epsilon$-symmetric Poincar\'e
nilcomplex if and only if $(S^{-1}D,\Phi)$ is an $\epsilon$-symmetric
Poincar\'e complex.  Similarly in the $\epsilon$-quadratic case.\\
{\rm (ii)} The nilpotent $\epsilon$-symmetric $L$-group of a ring with
involution $R$ fits into a split exact sequence:
$$0\to L_{\widetilde{K}_0(R)}^n(R[x],\epsilon) \to L^n_{\widetilde{K}_0(R)}
(R[x,x^{-1}],\epsilon) \to L\hbox{\rm Nil}^n(R,\epsilon) \to 0$$
with the surjection split by the injection
$$\begin{array}{l}
L\hbox{\rm Nil}^n(R,\epsilon) \to
L^n_{\widetilde{K}_0(R)}(R[x,x^{-1}],\epsilon)~;\\[1ex]
(C,\nu,\delta\phi,\phi) \mapsto (C[x,x^{-1}],[\nu,\delta\phi,\phi])
\oplus (C[x,x^{-1}],-\phi)\\[1ex]
([\nu,\delta\phi,\phi]_s=(x-\nu)\phi_s+T_{\epsilon}\delta\phi_{s-1}\,,\,s
\geqslant 0\,,\,\delta\phi_{-1}=0)~.
\end{array}$$
Similarly in the $\epsilon$-quadratic case, with a split exact sequence:
$$0\to L^{\widetilde{K}_0(R)}_n(R[x],\epsilon) \to L_n^{\widetilde{K}_0(R)}
(R[x,x^{-1}],\epsilon) \to L\hbox{\rm Nil}_n(R,\epsilon) \to 0$$
where the surjection split by the injection
$$\begin{array}{l}
L\hbox{\rm Nil}_n(R,\epsilon) \to
L^{\widetilde{K}_0(R)}_n(R[x,x^{-1}],\epsilon)~;\\[1ex]
(C,\nu,\delta\psi,\psi) \mapsto (C[x,x^{-1}],[\nu,\delta\psi,\psi])
\oplus (C[x,x^{-1}],-\psi)\\[1ex]
([\nu,\delta\psi,\psi]_s=
(x-\nu)\psi_s+T_{\epsilon}\delta\psi_{s+1}\,,\,s\geqslant 0)~.
\end{array}$$
{\rm (iii)} The morphism
$$\begin{array}{l}
L\hbox{\rm Nil}^n(R,\epsilon) \to L_{\widetilde{K}_0(R)}^n(R[x],\epsilon)~;~
(C,\nu,\delta\phi,\phi) \mapsto (C[x],\widetilde{\Phi})\\[1ex]
\hspace*{10mm}
(\widetilde{\Phi}_s=(1-x\nu)\phi_s+xT_{\epsilon}\delta\phi_{s-1}~,~s\geqslant
0~,~\delta\phi_{-1}=0)
\end{array}$$
is an isomorphism, and
$$L^n(R[x],\epsilon)~=~L^n(R,\epsilon)\oplus L\widetilde{\hbox{\rm
Nil}}^n(R,\epsilon)~,~
NL^n(R,\epsilon)~=~L\widetilde{\hbox{\rm Nil}}^n(R,\epsilon)~.$$
Similarly in the $\epsilon$-quadratic case, with the
morphism\footnote{As noted by the referee the cycles
$\widetilde{\Phi}, (1+T)\widetilde{\Psi} \in (W^{\%}C[x])_n$
differ by a boundary involving $\delta\psi_0$.}
$$\begin{array}{l}
L\hbox{\rm Nil}_n(R,\epsilon) \to L^{\widetilde{K}_0(R)}_n(R[x],\epsilon)~;~
(C,\nu,\delta\psi,\psi) \mapsto (C[x],\widetilde{\Psi})\\[1ex]
\hspace*{10mm}(\widetilde{\Psi}_s=
(1-x\nu)\psi_s+xT_{\epsilon}\delta\psi_{s+1}~,~s\geqslant 0)
\end{array}$$
an isomorphism, and
$$L_n(R[x],\epsilon)~=~L_n(R,\epsilon)\oplus L\widetilde{\hbox{\rm
Nil}}_n(R,\epsilon)~,~
NL_n(R,\epsilon)~=~L\widetilde{\hbox{\rm Nil}}_n(R,\epsilon)~.$$
\end{prop}
\noindent{\it Proof.} (i) Ranicki \cite{RHK}, Propositions 34.5.\\
(ii) The $\epsilon$-symmetric $L$-theory localization exact
sequence of Proposition 3.7.2 of Ranicki \cite{RES}
$$\dots  \to L_I^n(A,\epsilon) \to L_{S^{-1}I}^n(S^{-1}A,\epsilon)
\overset{\partial}{\to} L^n_I(A,S,\epsilon)
\to L^{n-1}_I(A,\epsilon)\to\dots$$
is defined for any ring with involution $A$, a central multiplicative
subset $S \subseteq A$ of nonzero divisors, and any $*$-invariant
subgroup $I \subseteq \widetilde{K}_0(A)$, with $L_I^n(A,S,\epsilon)$ the
cobordism group of $(n-1)$-dimensional $\epsilon$-symmetric Poincar\'e
complexes $(C,\phi)$ over $A$ such that
$$S^{-1}A\otimes_AC~ \simeq~ 0~,~[C]\in I~.$$
The boundary map is defined by
$$\partial~:~L_{S^{-1}I}^n(S^{-1}A,\epsilon) \to
L^n_I(A,S,\epsilon)~;~S^{-1}(C,\phi) \mapsto \partial (C,\phi)$$
with $(C,\phi)$ an $n$-dimensional $S^{-1}A$-Poincar\'e
$\epsilon$-symmetric complex over $A$ such that $[C] \in I$, and
$\partial (C,\phi)=(\partial C,\partial \phi)$ the
$(n-1)$-dimensional $S^{-1}A$-contractible $\epsilon$-symmetric
Poincar\'e complex over $A$ given by the boundary construction of
page 48 of \cite{RES}, with $\partial C={\mathcal
C}(\phi_0:C^{n-*} \to C)_{*+1}$.  For
$$(A,S)~=~(R[x],\{x^k\vert k \geqslant 0\})~~,~~S^{-1}A~=~R[x,x^{-1}]~,~
I~=~\widetilde{K}_0(R) \subseteq \widetilde{K}_0(R[x])$$
the localization exact sequence breaks up into split exact sequences
$$0  \to L^n_I(A,\epsilon) \to L^n_{S^{-1}I}(S^{-1}A,\epsilon)
\overset{\partial}{\longrightarrow}L^n_I(A,S,\epsilon) \to 0$$
with
$$\begin{array}{l}
L_I^n(A,\epsilon)~=~ L^n_{\widetilde{K}_0(R)}(R[x],\epsilon)~,\\[1ex]
L_{S^{-1}I}^n(S^{-1}A,\epsilon)~=~
L^n_{\widetilde{K}_0(R)}(R[x,x^{-1}],\epsilon)~,\\[1ex]
L_I^n(A,S,\epsilon)~=~L\hbox{\rm Nil}^n(R,\epsilon)
\end{array}$$
(Propositions 5.1.3, 5.1.4 of \cite{RES}). The formulae for
$[\nu,\delta\phi,\phi]$ and $[\nu,\delta\psi,\psi]$ are from page
445 of \cite{RES}. The identification
$L_I^n(A,S,\epsilon)=L\hbox{\rm Nil}^n(R,\epsilon)$ can be deduced
from (i), noting that by Proposition \ref{2.13} a finite f.g.
projective $R[x]$-module chain complex $D$ with projective class
$[D] \in I$ is such that $R[x,x^{-1}]\otimes_{R[x]}D \simeq 0$ if
and only if $D$ is chain equivalent to ${\mathcal C}(x-\nu:C[x]
\to C[x])$ for an $R$-nilcomplex $(C,\nu)$, with $C$ $R$-module
chain equivalent to
$D$ and $\nu \simeq x:C \simeq D \to C \simeq D$. The map
$$L{\rm Nil}^n(R,\epsilon) \to
L_I^n(R[x],S,\epsilon)~;~(C,\nu,\delta\phi,\phi)\mapsto (S^{-1}D,\Phi)$$
is an isomorphism, which factors as
$$L{\rm Nil}^n(R,\epsilon) \to L_{\widetilde{K}_0(R)}^n(R[x,x^{-1}],\epsilon)
\overset{\partial}{\longrightarrow} L_I^n(R[x],S,\epsilon)$$
with
$$\begin{array}{l}
L{\rm Nil}^n(R,\epsilon) \to
L_{\widetilde{K}_0(R)}^n(R[x,x^{-1}],\epsilon)~;\\[1ex]
(C,\nu,\delta\phi,\phi) \mapsto (C[x,x^{-1}],
\{(x-\nu)\phi_s+T_{\epsilon}\delta\phi_{s-1}\vert s\geqslant
0\})~~(\delta\phi_{-1}=0)~.
\end{array}$$
(iii)  The inclusion $R[x^{-1}] \to R[x,x^{-1}]$ induces a split injection
$$L^n_{\widetilde{K}_0(R)}(R[x^{-1}],\epsilon) \to L^n_{\widetilde{K}_0(R)}
(R[x,x^{-1}],\epsilon)~=~L^n_{\widetilde{K}_0(R)}(R[x],\epsilon)
\oplus L\hbox{\rm Nil}^n(R,\epsilon)$$
with image
$$L_p^n(R,\epsilon) \oplus L\widetilde{\hbox{\rm
Nil}}^n(R,\epsilon)~=~L{\rm Nil}^n(R,\epsilon)~.$$
Replacing $R[x^{-1}]$ by $R[x]$, it follows that the morphism
$$\begin{array}{l}
L{\rm Nil}^n(R,\epsilon) \to L^n_{\widetilde{K}_0(R)}
(R[x^{-1}],\epsilon)~;\\[1ex]
(C,\nu,\delta\phi,\phi) \mapsto (C[x^{-1}],
\{(1-x^{-1}\nu)\phi_s+x^{-1}T_{\epsilon}\delta\phi_{s-1}\vert
s\geqslant 0\})~~(\delta\phi_{-1}=0)
\end{array}$$
is an isomorphism. The inclusion $R[x^{-1}] \to R[x,x^{-1}]$
induces a split injection
$$L^n_{\widetilde{K}_0(R)}(R[x^{-1}],\epsilon) \to
L^n_{\widetilde{K}_0(R)}(R[x,x^{-1}],\epsilon)~=~
L^n_{\widetilde{K}_0(R)}(R[x],\epsilon)
\oplus L\hbox{\rm Nil}^n(R,\epsilon)$$
with image
$$L^n(R,\epsilon) \oplus L\widetilde{\hbox{\rm Nil}}^n(R,\epsilon)~,$$
and an isomorphism
$$L^n(R,\epsilon) \oplus L\widetilde{\hbox{\rm Nil}}^n(R,\epsilon) \to
L^n_{\widetilde{K}_0(R)}(R[x^{-1}],\epsilon)~.$$
\hfill\qed

In the applications of the nilpotent $L$-groups to the unitary
nilpotent $L$-groups we shall be particularly concerned with
the Witt groups of `nilforms' over $R$.

Define the {\it $Q{\rm Nil}$-groups} of an $R$-nilmodule $(P,\nu)$ to
be the $Q{\rm Nil}$-groups of the 0-dimensional $R$-nilcomplex
$(C,\nu^*)$ with
$$C~:~ \dots \to 0 \to C_0=P^* \to 0 \to \dots~,$$
as given in the $\epsilon$-symmetric case by
$$\begin{array}{ll}
Q{\rm Nil}^{\epsilon}(P,\nu)&=~Q{\rm Nil}^0(C,\nu^*,\epsilon)\\[1ex]
&=~\{\phi \in {\rm Hom}_R(P,P^*)\,\vert\, \epsilon \phi^*=\phi,
\nu^*\phi=\phi \nu:P \to P^*\}
\end{array}$$
and in the $\epsilon$-quadratic case by
$$\begin{array}{l}
Q{\rm Nil}_{\epsilon}(P,\nu)~=~Q{\rm Nil}_0(C,\nu^*,\epsilon)\\[1ex]
=~\displaystyle{\{(\delta\psi,\psi) \in {\rm Hom}_R(P,P^*)\oplus
{\rm Hom}_R(P,P^*)\,\vert\, \nu^*\psi-\psi \nu=\delta\psi+\epsilon
\delta\psi^*:P \to P^*\} \over
\{(\delta\chi-\epsilon\delta\chi^*+\nu^*\chi-\chi\nu,\chi-\epsilon\chi^*)
\,\vert\, (\delta\chi,\chi) \in {\rm Hom}_R(P,P^*)\oplus {\rm
Hom}_R(P,P^*)\}}~.
\end{array}$$
There is an evident $\epsilon$-symmetrization map
$$N_{\epsilon}~:~Q{\rm Nil}_{\epsilon}(P,\nu) \to Q{\rm
Nil}^{\epsilon}(P,\nu)~;~ (\delta\psi,\psi) \mapsto
N_{\epsilon}(\psi)~.$$

\begin{defn}\label{2.12} {\rm (\cite{RES}, p.452)\\
(i) A {\it nonsingular $\epsilon$-symmetric nilform over $R$}
$(P,\nu,\phi) $ consists of
\begin{itemize}
\item[(a)] an $R$-nilmodule $(P,\nu)$,
\item[(b)] an element $\phi \in Q{\rm Nil}^\epsilon(P,\nu)$ such that
$\phi:P \to P^*$ is an isomorphism.
\end{itemize}
Thus $(P,\phi)$ is a nonsingular $\epsilon$-symmetric form over $R$,
and there is defined an isomorphism of $R$-nilmodules
$$\phi~:~(P,\nu) \to (P^*,\nu^*)~.$$
A {\it lagrangian} for $(P,\nu,\phi)$ is a direct summand $L \subseteq
P$ such that
\begin{itemize}
\item[(c)] $\nu(L)\subseteq L$,
\item[(d)] the sequence
$$0 \to L\overset{i}{\to}P\overset{i^*\phi}{\longrightarrow}L^* \to 0$$
is exact, with $i:L \to P$ the inclusion.
\end{itemize}
In particular, $L$ is a lagrangian for the
nonsingular $\epsilon$-symmetric form $(P,\phi)$.\\
(ii) A {\it nonsingular $\epsilon$-quadratic nilform over $R$}
$(P,\nu,\delta\psi,\psi) $ consists of
\begin{itemize}
\item[(a)] an $R$-nilmodule $(P,\nu)$
\item[(b)] an element $(\delta\psi,\psi) \in
Q{\rm Nil}_{\epsilon}(P,\nu)$ such that
$N_\epsilon(\psi):P \to P^*$ is an isomorphism.
\end{itemize}
Thus $(P,\psi)$ is a nonsingular $\epsilon$-quadratic form over $R$,
and there is defined an isomorphism of $R$-nilmodules
$$N_{\epsilon}(\psi)~:~(P,\nu) \to (P^*,\nu^*)~.$$
A {\it lagrangian} for $(P,\nu,\delta \psi, \psi) $ is a direct summand
$L \subseteq P$ such that
\begin{itemize}
\item[(c)] $\nu(L)\subseteq L$,
\item[(d)] the sequence
$$0 \to
L\overset{i}{\to}P\overset{i^*N_\epsilon(\psi)}{\longrightarrow}L^* \to
0$$
is exact, with $i:L \to P$ the inclusion,
\item[(e)] $(i^*\delta\psi i,i^* \psi i)=(0,0) \in Q{\rm
Nil}_{\epsilon}(L,\nu\vert)$~.
\end{itemize}
In particular, $L$ is a lagrangian for the
nonsingular $\epsilon$-quadratic form $(P,\psi)$.\hfill$\qed$}
\end{defn}

The notion of stable isometry of nilforms is now defined in the usual
way using lagrangians and orthogonal direct sums, and $L\hbox{\rm
Nil}^0(R,\epsilon)$ (resp.  $L\hbox{\rm Nil}_0(R,\epsilon)$) is the
Witt group of nonsingular $\epsilon$-symmetric (resp.
$\epsilon$-quadratic) nilforms over $R$.  See Ranicki \cite{RES} (pp.
456-457) for the identification of $L\hbox{\rm Nil}^1(R,\epsilon)$
(resp.  $L\hbox{\rm Nil}_1(R,\epsilon)$) with the Witt group of
nonsingular $\epsilon$-symmetric (resp.  $\epsilon$-quadratic)
nilformations over $R$.

\begin{prop}\label{2.1}
{\rm (Ranicki \cite{RHK}, Proposition 41.3)}\\
{\rm (i)} For any ring with involution $R$ the
skew-suspension maps in the nilpotent $\pm\epsilon$-quadratic $L$-groups
are isomorphisms, so that
$$L\hbox{\rm Nil}_n(R,\epsilon)~=~L\hbox{\rm
Nil}_{n+2}(R,-\epsilon)~=~L\hbox{\rm Nil}_{n+4}(R,\epsilon)~,$$
with $L\hbox{\rm Nil}_{2n}(R,\epsilon)=L\hbox{\rm
Nil}_0(R,(-1)^n\epsilon)$ the Witt group of
nonsingular $(-1)^n\epsilon$-quadratic nilforms over $R$.
Similarly for $L\widetilde{\hbox{\rm Nil}}_*(R,\epsilon)$.\\
{\rm (ii)} If $R$ is a Dedekind ring with involution then
$$\begin{array}{l}
L\hbox{\rm Nil}^n(R,\epsilon)~=~L\hbox{\rm
Nil}^{n+2}(R,-\epsilon)~=~L\hbox{\rm Nil}^{n+4}(R,\epsilon)~,\\[1ex]
L\hbox{\rm
Nil}^n(R,\epsilon)~=~L_p^n(R,\epsilon)~~,~~L\widetilde{\hbox{\rm
Nil}}^n(R,\epsilon)~=~0~~(n \geqslant 0)~.
\end{array}$$
\end{prop}
\noindent{\it Proof.} (i) In order to establish the 4-periodicity
use algebraic surgery below the middle dimension, as for the
ordinary $\epsilon$-quadratic $L$-groups $L_n(R,\epsilon)$ in
Proposition I.4.3 of \cite{RATI,RATII} (cf.  Proposition
\ref{below} above).\\
(ii) The explicit proof in the case $n=0$ (\cite{RHK}, p.\ 588) extends
to the general case as follows.  Let $(C,\nu,\delta\phi,\phi)$ be an
$n$-dimensional $\epsilon$-symmetric Poincar\'e nilcomplex over $R$,
representing an element of $L\hbox{\rm Nil}^n(R,\epsilon)$, with
$$\nu^N~=~0~:~C \to C$$
for some $N \geqslant 1$.  We reduce to the case $N=1$ using the structure
theory of f.g.  modules over the Dedekind ring $R$ : every f.g.
$R$-module $M$ fits into a split exact sequence
$$0 \to T(M) \to M \to M/T(M) \to 0$$
with
$$T(M)~=~\{x \in M\,\vert\,ax =0 \in M~\hbox{for some}~a\neq 0 \in R\}$$
the torsion $R$-submodule and the quotient torsion-free $R$-module
$M/T(M)$ is f.g. projective. In particular, for any $R$-nilmodule $(P,\nu)$
with
$$\nu^N~=~0~:~P \to P$$
the $R$-submodule of $P$ defined by
$$T_N(P,\nu)~=~\{x \in P\,\vert\,ax \in \nu^{N-1}(P)~\hbox{for
some}~a\neq 0 \in R\}$$
is such that
$$T_N(P,\nu)/\nu^{N-1}(P)~=~T(P/\nu^{N-1}(P))~.$$
The torsion-free quotient $R$-module
$$(P/\nu^{N-1}(P))/T(P/\nu^{N-1}(P))~=~P/T_N(P,\nu)$$
is f.g. projective, so that $T_N(P,\nu)$ is a direct summand of $P$.
The inclusion defines a morphism of $R$-nilmodules
$$i~:~(T_N(P,\nu),0) \to (P,\nu)$$
Moreover, if $(P',\nu')$ is another $R$-nilmodule with $\nu'^N=0$ and
$$\theta~:~(P,\nu) \to (P',\nu')^*~=~({P'}^*,\nu'^*)$$
is a morphism of $R$-nilmodules then
$${i'}^*\theta i~=~0~:~T_N(P,\nu) \to T_N(P',\nu')^*$$
since for any $x\in T_N(P,\nu)$, $x' \in T_N(P',\nu')$ there exist
$a,a' \neq 0 \in R$, $y \in P$, $y' \in P'$ with
$$ax~=~\nu^{N-1}(y) \in P~~,~~a'x'~=~\nu'^{N-1}(x') \in P'$$
and
$$\begin{array}{ll}
a'\theta(x)(x')\overline{a}&=~\theta(ax)(a'x')\\[1ex]
&=~\theta(\nu^{N-1}(y))({\nu'}^{N-1}(y'))\\[1ex]
&=~\theta(\nu^{2N-2}(y))(y')\\[1ex]
&=~0 \in R~~\hbox{(since $2N-2 \geqslant N$)}
\end{array}$$
so that
$$\theta(x)(x')~=~0 \in R~.$$
Returning to the $n$-dimensional $\epsilon$-symmetric Poincar\'e
nilcomplex $(C,\nu,\delta\phi,\phi)$ with $\nu^N=0:C \to C$,
let  $i:(B,0) \to (C^{n-*},\nu^*)$ be the inclusion of the subcomplex
defined by
$$B_r~=~T_N(C^{n-r},\nu^*)~.$$
The chain map of $R$-nilmodule chain complexes defined by
$$f~=~i^*~:~ (C,\nu) \to (D,0)~=~(B^{n-*},0)$$
is such that
$$f^*(\delta\phi,\phi)~=~0 \in Q\hbox{\rm Nil}^n(D,0,\epsilon)~.$$
Algebraic surgery on $(C,\nu,\delta\phi,\phi)$ using the $(n+1)$-dimensional
$\epsilon$-symmetric nilpair $(f:(C,\nu) \to (D,0),(0,(\delta\phi,\phi)))$
over $R$ results in a cobordant
$n$-dimensional $\epsilon$-symmetric Poincar\'e nilcomplex
$(C',\nu',\delta\phi',\phi')$ over $R$ with
$$\nu'~ \simeq~ 0~:~C' \to C'~.$$
\hfill\qed

\subsection{ The unitary nilpotent $L$-groups $\hbox{\rm UNil}$}

Let $R$ be any ring.
An {\it involution} on an $R$-$R$ bimodule $\cA$ is a homomorphism
$$\cA\to \cA~;~a\mapsto \bar a$$
which satisfies
$$\bar {\bar a}~=~a~,~\overline{ras}~=~\bar s \bar a \bar r~\hbox{\rm
for all } a\in \cA,\; r,s\in R~.$$
For any $R$-module $P$ there is defined an $R$-module
$$\cA P~=~\cA\otimes_RP~.$$
As in the special case $\cA=R$ write the evaluation pairing as
$$\langle~,~\rangle~:~ \cA P^* \times P  \to \cA~;~
(a\otimes f, x) \mapsto \langle a \otimes f,x \rangle~=~af(x)~.$$
An element $\phi \in \hbox{\rm Hom}_R(P,\cA P^*) $ determines a
$\cA$-valued sesquilinear form on $P $
$$\langle~,~\rangle_\phi~:~P \times P \to \cA~;~(x,y) \mapsto \langle
\phi(x),y \rangle~,$$
and we identify $\hbox{\rm Hom}_R(P,\cA P^*) $ with the additive group
of such forms.  For $\epsilon=\pm 1$ and a f.g.  projective $P$ define
an involution
$$T_{\epsilon}~:~\hbox{\rm Hom}_R (P,\cA P^*)\to \hbox{\rm Hom}_R
(P,\cA P^*)~;~\phi \mapsto \epsilon \phi^*~~,~~
\langle x,y \rangle_{\phi^*}~=~\overline{\langle y,x \rangle}_\phi~.$$
One then defines a map
\begin{equation}\label {2.2}
N_{\epsilon}~=~1+T_{\epsilon}~:~\hbox{\rm Hom}_R (P,\cA P^*)\to
\hbox{\rm Hom}_R (P,\cA P^*)~;~\phi \mapsto \phi +\epsilon \phi^t
\end{equation}
with
$$\langle x,y \rangle_{N_{\epsilon}(\phi)}~=~
\langle x,y \rangle_{\phi}+\epsilon \overline{\langle y,x \rangle}_{\phi}~.$$
An $\cA$-valued $\epsilon$-symmetric form $(P,\lambda)$
(resp. $\epsilon$-quadratic form $(P,\mu)$) over $R$
is a f.g. projective $R$-module $P$ together with an element of the group
$$\begin{array}{l}\label{2.3}
\lambda \in Q^\epsilon (P,\cA) ~=~
\hbox{\rm ker}(1-T_{\epsilon}:
\hbox{\rm Hom}_R (P,\cA P^*)\to \hbox{\rm Hom}_R (P,\cA P^*))~,\\[1ex]
\mu \in Q_\epsilon(P,\cA) ~=~
\hbox{\rm coker}(1-T_{\epsilon}:
\hbox{\rm Hom}_R (P,\cA P^*)\to \hbox{\rm Hom}_R (P,\cA P^*))~.
\end{array}$$
As usual, for $\lambda \in Q^\epsilon(P,\cA)$ we write
$$\lambda(x,y)~=~\langle \lambda(x),y \rangle \in \cA$$
and for $\mu \in Q_\epsilon (P,\cA) $ we write
$$\mu(x)~=~
\langle \mu(x),x \rangle \in  \cA/\{a-\epsilon \bar{a}|a\in \cA\}~.$$
The map $N_\epsilon $ induces a well defined map:
\begin{equation}\label {2.5}
N_{\epsilon}~:~Q_\epsilon(P,\cA) \to Q^\epsilon(P,\cA)~;~
[\mu]\mapsto \mu +\epsilon \mu^t~.
\end{equation}

\begin{defn}\label{2.6}{\rm (Cappell \cite{CaB})\\
(i) Let $\cB_1,\cB_{-1} $ be $R$-bimodules with involution.
Assume $\cB_1,\cB_{-1} $ are free as right $R$-modules.
A {\it nonsingular $\epsilon$-quadratic unilform over} $(R;\cB_1,\cB_{-1})$
is a quadruple
$$(P_1, P_{-1},\mu_1,\mu_{-1})$$
where, for $\delta = \pm 1$,  we require:
\begin{itemize}
\item[(a)] $(P_\delta,\mu_\delta) $ is a stably f.g. free $\cB_\delta$-valued
$\epsilon$-quadratic form over $R$,
\item[(b)] $P_{\delta}=P_{-\delta}^* $; we then identify
$(P_{\delta}^*)^*=P_\delta$
in the usual way, and write the evaluation pairing as
$$\langle~,~\rangle~:~P_1\times P_{-1}\to R ~;~(x,f) \mapsto f(x)~.$$
\item[(c)] If $\lambda_\delta = N_\epsilon (\mu_\delta)$ is the
associated $\epsilon$-symmetric form to $\mu_\delta$, then the
composite
$$P_1\overset{\lambda_1}{\to}
\cB_1P_{-1}\overset{\lambda_{-1}\otimes1}{\longrightarrow}
\cB_{-1}\cB_1P_{1}\overset{\lambda_{1}\otimes1}{\longrightarrow}
\cB_1\cB_{-1}\cB_1P_{-1}
\overset{\lambda_{-1}\otimes1}{\longrightarrow}\dots$$
is eventually zero. (That is to say, for some $k$, the composite map
$P_1 \to (\cB_{-1}\cB_1)^kP_1$ is zero.)
\end{itemize}
(ii) A {\it sublagrangian} for $(P_1,P_{-1},\mu_1,\mu_{-1})$
is a pair of stably f.g. free direct summands $V_1\subseteq P_1$,
$V_{-1}\subseteq P_{-1}$ such that, for $\delta = \pm 1$
\begin{equation}\label{2.7}
\langle V_1, V_{-1}\rangle ~=~ 0~~,~\lambda_\delta(V_\delta)
\subseteq \cB_\delta V_{-\delta}~~,~~\mu_\delta(x)~=~0~ \hbox{\rm
for~ all}~x \in V_\delta~.
\end{equation}
We call $(V_1, V_{-1}) $ a {\it lagrangian} if in addition:
\begin{equation}\label{2.8}
V_1 ~=~ V_{-1}^{\perp}.
\end{equation}
\hfill$\qed$}
\end{defn}

One can form orthogonal direct sums of $\epsilon$-quadratic
unilforms over $(R;\cB_1,\cB_{-1})$ in a rather obvious way.
Cappell \cite{CaB} defined $\hbox{\rm
UNil}_{2n}(R;\cB_1,\cB_{-1})$ to be the Witt group of stable
isometry classes of nonsingular $(-1)^n$-quadratic unilforms over
$(R;\cB_1,\cB_{-1})$ modulo those admitting lagrangians, and
showed (geometrically) that if $\pi_{-1},\pi_0,\pi_1$ are finitely
presented groups with $\pi_0 \subseteq \pi_{-1}$, $\pi_0 \subseteq
\pi_1$ and
$$\pi~=~\pi_{-1}*_{\pi_0}\pi_1~~,~~
R~=~\Z[\pi_0]~,~{\cB}_{\pm 1}~=~\Z[\pi_{\pm 1}\backslash \pi_0]$$
then the morphism defined by
$$\begin{array}{l}
\hbox{\rm UNil}_{2n}(R;\cB_1,\cB_{-1}) \to L_{2n}(\Z[\pi])~;\\[1ex]
(P_1, P_{-1},\mu_1,\mu_{-1}) \mapsto
(\Z[\pi]\otimes_{\Z[\pi_0]}(P_1\oplus P_{-1}),
\begin{pmatrix}
\mu_1 & 1 \\
0 &  \mu_{-1}
\end{pmatrix})
\end{array}$$
is a split monomorphism.

If an $\epsilon$-quadratic unilform $u=(P_1,P_{-1},\mu_1,\mu_{-1})$ has
a sublagrangian $(V_1,V_{-1})$, then one can form a new $\epsilon$-quadratic
unilform (see Connolly and Ko\'zniewski \cite{CK}, 6.3 (f))
$$u'~=~(V_{-1}^\perp/V_1, V^\perp_{1}/V_{-1},\mu'_1,\mu'_{-1})~,$$
so that
$$[u] ~=~ [u'] \in \hbox{\rm UNil}_{2n}(R;\cB_1,\cB_{-1})~.$$

\subsection{The proof of Theorem A in the even-dimensional case.}

We begin by defining maps:
$$
L\widetilde{\hbox{\rm Nil}}_{2n}(R) \overset{c}{\to}\hbox{\rm
UNil}_{2n}(R;R,R)\overset{r}{\to} NL_{2n}(R)\subseteq L_{2n}(R[x]).
$$
The proof will show that the maps $c,r$ are both isomorphisms.

Let $\epsilon = (-1)^n.$
\begin{defn}\label{2.15}
{\rm
The map
$$r~:~\hbox{\rm UNil}_{2n}(R;R,R) \to NL_{2n}(R)~;~u \mapsto r(u)$$
sends an $\epsilon$-quadratic unilform  $u=(P_1,P_{-1},\mu_1,\mu_{-1})$
over  $(R;R,R)$ to the $\epsilon$-quadratic form $r(u)$ over $R[x]$ given by:
$$r(u) ~=~ (P_1[x]\oplus P_{-1}[x],\psi_0 +x\psi_1)$$
where
$$\psi_0 ~=~\begin{pmatrix}
0&1\\
0&\mu_{-1}
\end{pmatrix}~,~
\psi_1 ~=~\begin{pmatrix}
\mu_1& 0\\
0&0
\end{pmatrix}~.$$}
\hfill$\qed$
\end{defn}
Here, $\psi_i:(P_1\oplus P_{-1})[x]\to
(P_{-1}\oplus P_1)[x] $ $(i=0,1)$
is the $R[x]$-module morphism induced, using change of
coefficients,  from the $R$-module morphism of the same name
$$\psi_i~:~(P_1\oplus P_{-1}) \to (P_1\oplus P_{-1})^*~=~(P_{-1}\oplus P_1)~.$$
In order to verify that $r$ is well-defined, first notice that
$$N_\epsilon(\psi_0+x\psi_1) ~=~
\begin{pmatrix}
0&1\\
\epsilon & 0
\end{pmatrix}
(1+\nu)~:~(P_1\oplus P_{-1})[x] \to (P^*_1\oplus P^*_{-1})[x]$$
where
$$\nu~=~ \left(
\begin{matrix}
0&\epsilon \lambda_{-1}\\
x\lambda_1 & 0
\end{matrix}
\right)~:~(P_1\oplus P_{-1})[x] \to (P_1\oplus P_{-1})[x]~
,\quad \lambda_{\pm 1}: ~=~ N_\epsilon(\mu_{\pm 1})~.$$
Because
$$\nu^2~=~ \left(
\begin{matrix}
x\epsilon \lambda_{-1}\lambda_1&0 \\
0&x\lambda_1\lambda_{-1}
\end{matrix}
\right)$$
Definition \ref{2.6} shows that $\nu$ is obviously
nilpotent. Therefore $N_\epsilon(\psi_0 +x\psi_1)$ is nonsingular.

To see that $[r(u)]\in NL_{2n}(R)$, notice that
$\eta_*[r(u)]=[P_1\oplus P_{-1},\psi_0]$, and that $P_1\oplus 0$ is a
lagrangian for  $(P_1\oplus P_{-1},\psi_0)$.

The rule $u\mapsto r(u)$ preserves orthogonal direct sums of
forms. If $(V_1,V_{-1}) $ is a lagrangian for $u$, then
$V_1[x]\oplus V_{-1}[x] $ is a lagrangian for $r(u)$.
We thus have a well-defined homomorphism:
$$r~:~\hbox{\rm UNil}_{2n}(R;R,R) \to NL_{2n}(R)~.$$

\begin{defn}\label{2.17} {\rm The map
$$c~:~ L\widetilde{\hbox{\rm Nil}}_{2n}(R) \to \hbox{\rm
UNil}_{2n}(R;R,R)~;~z \mapsto c(z)$$
sends a nonsingular
$\epsilon$-quadratic nilform $z=(P,\nu,\delta\psi,\psi)$ over $R$
(see Definition \ref{2.12}) to $c(z)= (P_1, P_{-1},\mu_1,\mu_{-1})$,
where
$$P_1~=~P~~,~~P_{-1}~=~P^*~~,~~\mu_1~=~\delta\psi -\nu^*\psi~~,~~
\mu_{-1}~=~-\phi^{-1}\psi^* \phi^{-1}$$
with $\phi=N_\epsilon(\psi):P\to P^*$ an isomorphism.\hfill$\qed$ }
\end{defn}
Using Definition \ref{2.12} set
$$\lambda_1~=~N_\epsilon(\mu_1) ~=~ -\phi\nu ~=~ -\nu^*\phi~,$$
noting that $N_{-\epsilon}N_\epsilon(\delta\psi)=0$. Set also
$$ \lambda_{-1} ~=~ N_\epsilon(\mu_{-1}) ~=~-\epsilon\phi^{-1}~.$$
Because $\lambda_{-1} \lambda_1 = \epsilon \nu$, and $\nu$ is
nilpotent, it follows that $c(z)$ is an $\epsilon$-quadratic unilform over
$(R;R,R).$ The rule $ z\mapsto c(z) $ preserves orthogonal
direct sums. Moreover, if $N$ is a lagrangian for $z $, then
$(N, N^\perp)$ is a lagrangian for $c(z)$. Therefore
Definition \ref{2.17} gives a homomorphism:
\begin{equation}\label{2.18}
c~:~L\widetilde{\hbox{\rm Nil}}_{2n}(R) \to \hbox{\rm UNil}_{2n}(R;R,R)
\end{equation}

\begin{defn}\label{2.19}
{\rm The morphism
$$j~:~L\widetilde{\hbox{\rm Nil}}_{2n}(R) \to NL_{2n}(R)~;~y \mapsto j(y)$$
sends $y=[P,\nu,\delta\psi,\psi]$ to
$$j(y) ~=~ [P[x],\psi +x(\delta\psi-\nu^*\psi)]~.\eqno{\qed}$$}
\end{defn}

It was proved in Ranicki \cite{RES}, p.\ 445 that $j$ is in fact an
isomorphism.  See Remark \ref{2.21} below for the precise matching up
of the formula in Definition \ref{2.19} with the morphism defined
there.

The right hand side in Definition \ref{2.19} gives a nonsingular form because:
$$
N_\epsilon(\psi+x(\delta\psi-\nu^*\psi)) ~=~ N_\epsilon(\psi)(1-x\nu),
$$
an isomorphism by Definition \ref{2.12}. Moreover this right
hand side is in $NL_{2n}(R)$, also by Definition \ref {2.12}.

\begin{rem}\label{2.21}
{\rm In order to obtain the formula in Definition \ref{2.19} for
$j(y)$ from the formula in \cite{RES}, p.\ 445 one must make the
following translation of the terminology there to our terminology:
$$\begin{array}{l}
A~=~R~~,~~C^0 ~=~ P~~,~~C^i ~=~0~\hbox{\rm for}~i\neq 0~,\\[1ex]
\psi_0 ~=~\psi~~,~~\delta\psi_1~=~ \delta \psi~,
\end{array}$$
noting that the $x^{-1}$ is our $x$, and the $\nu^*$ there is our
$\nu$.  In the following argument we shall use the Witt group
 $$L\hbox{\rm Nil}^h_{2n}(R)~=~L_{2n}(R) \oplus
 L\widetilde{\hbox{\rm Nil}}_{2n}(R)$$
of nonsingular $(-)^n$-quadratic nilforms
$(P,\nu,\delta\psi,\psi)$ over $R$ with $P$ a f.g. free
$R$-module, and the split injection
$$\Delta~:~ L\hbox{\rm Nil}^h_{2n}(R) \to L_{2n}(R[x,x^{-1}])~;~
[P,\nu,\delta\psi,\psi] \mapsto [P[x,x^{-1}],(x^{-1}-\nu^*)\psi+\delta\psi]$$
defined there, along with the splitting map
$$\partial~:~L_{2n}(R[x,x^{-1}]) \to L\hbox{\rm Nil}^h_{2n}(R)$$
and the natural inclusion and projection:
$$L\widetilde{\hbox{\rm Nil}}_{2n}(R)\overset{i}{\to}
L\hbox{\rm Nil}^h_{2n}(R)\overset{p}{\to}L\widetilde{\hbox{\rm
Nil}}_{2n}(R)~.$$
\indent Let $\tilde{E}:NL_{2n}(R)\to
L_{2n}(R[x,x^{-1}]) $ be the restriction of the natural
monomorphism
$$E~:~L_{2n}(R[x])\to L_{2n}(R[x,x^{-1}])~.$$
Also, set
$$\begin{array}{l}
\tilde\partial~=~p\partial~:~L_{2n}(R[x,x^{-1}]) \to
L\widetilde{\hbox{\rm Nil}}_{2n}(R)~,\\[1ex]
\tilde\Delta~=~\Delta i~:~ L\widetilde{\hbox{\rm Nil}}_{2n}(R) \to
L_{2n}(R[x,x^{-1}])~.
\end{array}$$
Because $\p\Delta = 1, $ we get $\tilde{\p}\tilde{\Delta} = 1$.
According to the braid on page 448 of \cite{RES},
$\tilde{\p}\tilde{E}$ is an isomorphism. The map $j$ of Definition
\ref{2.19} is $j=(\tilde{\p}\tilde{E})^{-1}$. To get the
formula for $j$ in Definition \ref{2.19}, note that the
''devissage'' map $\tilde{\p}$ satisfies:
$$\tilde{\p}~=~\tilde{\p}M$$
with
$$M~:~L_n(R[x,x^{-1}])\to
L_n(R[x,x^{-1}])~;~(P,\psi) \mapsto (P,x\psi)~.$$
Then from \cite{RES}, p. 445, we translate and find:
$$\begin{array}{ll}
\tilde{\Delta} (y) &=~ \tilde{\Delta}([P,\nu,\delta\psi,\psi])\\[1ex]
&=~ [P[x,x^{-1}], x^{-1}\{\psi +x(\delta\psi-\nu^*\psi)\}]~.
\end{array}$$
So
$$\begin{array}{ll}
j(y)&=~ j(\tilde{\p}M\tilde{\Delta}(y))\\[1ex]
&=~\tilde{E}^{-1}M\tilde{\Delta}(y)\\[1ex]
&=~\tilde{E}^{-1}([P[x,x^{-1}],\psi+x(\delta\psi-\nu^*\psi)])\\[1ex]
&=~[P[x],\psi+x(\delta\psi-\nu^*\psi)]~,
\end{array}$$
as in Definition \ref{2.19}.}\hfill$\qed$\end{rem}

As explained above, \cite{RES} proves that $j$ is an isomorphism.

\begin{rem}{\rm The inverse of $j$
\begin{equation} \label{2.22}
k~=~j^{-1}~:~NL_{2n}(R) \to L\widetilde{\hbox{\rm Nil}}_{2n}(R)\end{equation}
can be computed via Higman linearization (see Connolly and Ko\'zniewski
\cite{CK}, 3.6 (a)) in the following way.  By Higman linearization,
each element of
$NL_{2n}(R)$ can be represented in the form $[P[x],
\psi_0+x\psi_1]$. In these terms, the formula for $k=j^{-1}$ is:
\begin{equation}\label{2.23}
k[P[x],\psi_0+x\psi_1] ~=~ [P,\nu,\delta\psi,\psi]~,
\end{equation}
where
$$\psi ~=~ \psi_0~~,~~\nu ~=~(N_\epsilon(\psi_0))^{-1}N_\epsilon(\psi_1)~~,~~
\delta \psi ~=~ \nu^*\psi_0 +\psi_1~.$$
It is clear that $jk=1$.}\hfill$\qed$
\end{rem}

We now turn to the  proof of Theorem A in even
dimensions. We only have to show that:
\begin{equation}\label{2.24}
{\rm (i)}~ckr ~=~ 1~~,~~{\rm (ii)}~~rc ~=~ j~.
\end{equation}
The proof of (\ref{2.24}) (i) is easiest:
let $(P_1,P_{-1}, \mu_1,\mu_{-1})$ be an $\epsilon$-quadratic unilform over
$(R;R,R)$.  By Definitions \ref{2.15}, \ref{2.17} and (\ref{2.23}), and
direct calculation, we obtain:
\begin{equation}\label{2.25}
ckr[P_1,P_{-1},\mu_1,\mu_{-1}]~=~[P_1\oplus P_{-1},
P_{-1}\oplus P_1,\tilde{\mu}_1,\tilde{\mu}_{-1}]
\end{equation}
where
$$\tilde{\mu}_1~=~
\begin{pmatrix}
\mu_1&0\\
    0&0
\end{pmatrix},
\tilde{\mu}_{-1}~=~
\begin{pmatrix}
\mu_{-1}&1\\
    0&0
\end{pmatrix}~.$$
Perform a sublagrangian construction on the right hand side of
(\ref{2.25}), using the sublagrangian
$$V_1~=~ 0\oplus P_{-1}~,~V_{-1} ~=~ 0~.$$
This yields:
$$
[P_1\oplus P_{-1}, P_{-1}\oplus P_1,\tilde{\mu}_1,
\tilde{\mu}_{-1}]~=~ [P_1,P_{-1},
\mu_1,\mu_{-1}].
$$
Therefore $ckr=1$, proving equation (\ref{2.24}) (i).

Next we prove equation (\ref{2.24}) (ii).

Suppose $a= [P,\nu,\delta\psi,\psi]\in
L\widetilde{\hbox{\rm Nil}}_{2n}(R) $. By direct calculation and
Definitions \ref{2.15}, \ref{2.17}, we have
\begin{equation}\label {2.26}
rc(a) ~=~ [P[x]\oplus P^*[x],\Psi_0+x\Psi_1]\end{equation}
where
$$\Psi_0~=~ \begin{pmatrix}
0&1\\
0&-\phi^{-1}\psi^*\phi^{-1}
\end{pmatrix}~~,~~
\Psi_1~=~
\begin{pmatrix}
\delta\psi-\nu^*\psi&0\\
0&0
\end{pmatrix}$$
with $\phi=N_\epsilon(\psi)$.  By hypothesis (see Definition
\ref{2.12}), $(P,\psi) $ admits a lagrangian, say $N\subseteq P$. Let
$$V~=~ (\phi N)[x]\subseteq P^*[x]\subseteq P[x]\oplus P^*[x]~.$$
By (\ref{2.26}) $V$ is a sublagrangian for $\Psi_0+x\Psi_1 $.
In fact, setting $\Phi= N_\epsilon(\Psi_0+x\Psi_1) $, one readily
computes that the $\Phi$-orthogonal complement of $V$ is
$$V^\perp_\Phi ~=~ \{(u,v)\in P[x]\oplus P^*[x]\,|\,\phi(u)-v\in V\}~.$$
Therefore one obtains an isomorphism
$$g~:~P[x] \to  V^\perp_\Phi/V~;~u \mapsto (u,\phi(u))~.$$
Let $ (V^\perp_\Phi/V,\Psi') $ be the sublagrangian
construction on $cr(a)$ using $V$.  We claim that
\begin{equation}\label{2.29}
g~:~(P[x],\psi +x(\delta\psi -\nu^*\psi))\to
(V^\perp_\Phi/V,\Psi')
\end{equation}
is an isometry.  Since the right hand side of (\ref{2.29}) represents
$rc(a)$, and the left hand side is $j(a)$, this claim (\ref{2.29}) will
prove (\ref{2.24}) (ii).

We prove (\ref{2.29}) using the duality pairing
$$\begin{array}{ll}
\{~,~\}~:&(P^*[x]\oplus P[x]) \times (P[x]\oplus P^*[x]) \to R[x]~;\\[1ex]
&((\xi,\eta), (\eta',\xi')) \mapsto
\{(\xi,\eta), (\eta',\xi')\}~=~
\langle\xi,\eta'\rangle+\overline{\langle\xi', \eta\rangle}~.
\end{array}$$
(\ref{2.29}) amounts to the identity:
\begin{equation}\label{2.30}
\langle [\psi +x(\delta\psi-\nu^*\psi)](u),v\rangle ~=~ \{\Psi(u,\phi(u)),
(v,\phi(v))\}~~(u,v\in P[x])~,
\end{equation}
where $\Psi=\Psi_0 +x\Psi_1:P[x]\oplus P^*[x] \to P^*[x]\oplus P[x]$ .
The right hand side of (\ref{2.30}) is computed from (\ref{2.26}) as:
$$\begin{array}{l}
\langle[\phi +x(\delta\psi-\nu^*\psi)](u),v\rangle
+\overline{\langle\phi(v),-\phi'\psi^*(u)\rangle}\\[1ex]
\hspace*{20mm}=~
\langle[\phi+x(\delta\psi-\nu^*\psi](u),
v\rangle+\langle-\phi^*\phi^{-1}\psi^*(u),v\rangle\\[1ex]
\hspace*{20mm}=~\langle[(\phi-\epsilon\psi^*)+
x(\delta\psi-\nu^*\psi)](u),v\rangle,
\end{array}$$
which is the left hand side of (\ref{2.30}). This proves
(\ref{2.29}) and therefore also (\ref{2.24}) (ii). Therefore
the proof of Theorem A, when $n$ is even, is complete.
\begin{rem} {\rm
It seems appropriate to record here an explicit formula for
the inverse isomorphism
$$c^{-1}~:~\hbox{\rm UNil}_{2n}(R;R,R) \to L\widetilde{\hbox{\rm
Nil}}_{2n}(R)$$
which can be derived from (\ref{2.23}), (\ref{2.24}) and
Definition \ref{2.15}, as follows.

For $ [P_1,P_{-1},\mu_1,\mu_{-1}]\in \hbox{\rm UNil}_{2n}(R;R,R)$ we have:
\begin{equation}\label{2.32}
c^{-1}([P_1,P_{-1},\mu_1,\mu_{-1}]) ~=~
(P_1\oplus P_{-1},\nu,\delta\psi,\psi)
\end{equation}
where
$$\begin{array}{l}
\psi ~=~\begin{pmatrix}
0&1\\
0&\mu_{-1}
\end{pmatrix}
:P_1\oplus P_{-1}\to P_{-1}\oplus P_1~~,~~\psi_1 ~=~
\begin{pmatrix}
\mu_1&0\\
0&0
\end{pmatrix}~,\\[1ex]
\delta\psi~=~\nu^*\psi +\psi_1~:~P_1\oplus P_{-1}\to P_{-1}\oplus P_1~,\\[1ex]
\nu~=~-N_\epsilon(\psi)^{-1}N_\epsilon(\psi_1) ~=~
\begin{pmatrix}
\epsilon\lambda_{-1}\lambda_1 & 0\\
-\lambda_1& 0
\end{pmatrix}
\end{array}$$
with $\lambda_{\pm 1}=N_\epsilon(\mu_{\pm 1})$.\hfill$\qed$}
\end{rem}

\subsection{The proof of Theorem A in the odd-dimensional case}

We begin by commenting that the ``simple $L$-theory'' version of
Theorem A, in even dimensions,  proceeds uneventfully, along
the same lines as above. We explain this in some detail now.

$\hbox{\rm UNil}^s_{2n}(R;\cB_1,\cB_{-1}) $ is defined in Cappell
\cite{CaB} (p.1118). Also,
$$L\hbox{\rm Nil}^s_n(R)~=~ L^s_n(R)\oplus L\widetilde{\hbox{\rm Nil}}^s_n(R)$$
is defined in (\cite{RES}, p. 466-468), where there are also
constructed exact sequences:
$$\begin{array}{l}
0\to L_n^{I_+}(R[x])\overset{E^s}{\longrightarrow}
L^J_n(R[x,x^{-1}])\longrightarrow
L^p_n(R)\oplus L\widetilde{\hbox{\rm Nil}}^s_n(R)\to 0\\[1ex]
0\to L_n^{I_-}(R[x^{-1}])\longrightarrow
L^J_n(R[x,x^{-1}])\overset{\p^s}{\longrightarrow}
L^p_n(R)\oplus L\widetilde{\hbox{\rm Nil}}^s_n(R)\to 0
\end{array}$$
where
$$\begin{array}{l}
I_{\pm}~=~\widetilde{K}_1(R)\subseteq \widetilde{K}_1(R[x^{\pm
1}])~=~ \widetilde{K}_1(R) \oplus \widetilde{\rm Nil}_0(R)~,\\[1ex]
J~=~\widetilde{K}_1(R)\oplus K_0(R) \subseteq
\widetilde{K}_1(R[x,x^{-1}]) ~=~\widetilde{K}_1(R)\oplus K_0(R)
\oplus \widetilde{\rm Nil}_0(R) \oplus \widetilde{\rm Nil}_0(R)~.
\end{array}$$
Define
$$NL^s_n(R)~=~\hbox{\rm ker}(\eta_*:L_n^{K_1(R)}(R[x]) \to L_n(R))$$
and let $\widetilde{\Delta}^s,\tilde{E}^s,\tilde{\p}^s $ be as
in Remark \ref{2.21}, concluding that $\tilde{\p}^s\tilde{E}^s$
is an isomorphism. As before, define
$$j^s =(\tilde{\p}^s\tilde{E}^s)^{-1}~:~
L\widetilde{\hbox{\rm Nil}}_{2n}(R)\to NL^s_{2n}(R)$$
using the formula of Definition \ref{2.19}. The maps
$$
L\widetilde{\hbox{\rm Nil}}^s_{2n}(R)\overset{c^s}{\longrightarrow}
\hbox{\rm UNil}_{2n}^s(R;R,R)\overset{r^s}{\longrightarrow}NL^s_{2n}(R),
\hbox{\rm  and } k^s~=~(j^s)^{-1}
$$
are now defined exactly as in Definitions \ref{2.15}--\ref{2.23}, and
the proof that these are isomorphisms can now be repeated
without change. In summary, we have:
\begin{prop}\label{2.34}
The maps  $L\widetilde{\hbox{\rm Nil}}^s_{2n}(R)\overset{c^s}{\longrightarrow}
\hbox{\rm UNil}_{2n}^s(R;R,R)\overset{r^s}{\longrightarrow}NL^s_{2n}(R)$
described in the paragraph above are isomorphisms. Moreover,
$j^s=r^sc^s.$
\end{prop}

We now complete the proof of Theorem A in odd dimensions.

Let $S=R[z,z^{-1}]$, extending the involution on $R$ to $S$ by
$$\overline{z}~=~ z^{-1}~.$$
Let $i:R \to S$ be the inclusion. The split exact
sequence of Shaneson \cite{S} and Ranicki \cite{RII}
$$0\to L^s_n(R)\overset{i_*}{\to}L^s_n(S)\to L_{n-1}(R) \to 0$$
yields the split exact sequence
$$0\to NL^s_n(R)\overset{i_*}{\to}NL^s_n(S)\to NL_{n-1}(R) \to 0~.$$
Cappell \cite{CaB} defined $\hbox{\rm UNil}_{2n-1}(R;R,R)$ as the cokernel
in the split exact sequence:
$$
0\to \hbox{\rm UNil}^s_{2n}(R;R,R)\to \hbox{\rm UNil}^s_{2n}(S;S,S)\to
\hbox{\rm UNil}_{2n-1}(R;R,R)\to 0.
$$
The isomorphism $r^s$ of Proposition \ref{2.34}, being
functorial, therefore induces an isomorphism:
$$r~:~\hbox{\rm UNil}_{2n-1}(R;R,R) \to NL_{2n-1}(R)~.$$
This proves Theorem A.

\section{Chain bundles and the proof of Theorem B.}\label{sect3}

\subsection{Universal chain bundles.}

We begin with a resum\'e of the results of Ranicki
\cite{RATI,RATII},\cite{RAP} and Weiss \cite{WeI,WeII} which we
need.  As in Section \ref{sect2}, $R$ is a ring with involution.

A {\it chain bundle} $(B,\beta)$ over $R$ is a projective $R$-module chain
complex $ B $ together with a 0-cycle
$$\beta\in(\widehat W^\%B^{-*})_0~.$$
(We shall be mainly concerned with cases when the chain modules $B_r$
are f.g. projective.)
A map of chain bundles $f:(C,\gamma)\to (B,\beta)$ is a chain map
$f:C\to B $ such that
$$[\widehat f^\%(\beta)] ~=~ [\gamma] \in \widehat{Q}^0(C^{-*})~,$$
with $\widehat f^\%:\widehat W^\%B^{-*}\to W^\%C^{-*}$ the chain map
induced by $f$.  Each chain bundle $(B,\beta)$ determines a homomorphism
\begin{equation}\label{3b}
J_{\beta}~:~Q^n(B) \to \widehat Q^n(B)~;~
\phi \mapsto J(\phi)-\widehat \phi_0^\%(S^n\beta)
\end{equation}
where $J$ is as in (\ref{3a}), $\phi_0 ^\%$ is the map induced by
$\phi_0:B^{n-*}\to B$, and  $S^n: \widehat W^\%C\to \Sigma^{-n}\widehat
W^\%\Sigma^n(C)$ is the natural isomorphism of chain complexes.
The map $J_\beta$ is not induced by a chain map.

The Tate $\Z_2$-cohomology group
$$\widehat H^r(\Z_2;R)~=~\{x \in R\,\vert\,\bar{x}=(-1)^rx \}/
    \{y +(-1)^r \bar{y}\,\vert\, y \in R\}$$
is an $R$-module via
$$R \times \widehat H^r(\Z_2;R) \to \widehat H^r(\Z_2;R)~;~
(a,x) \mapsto a x \bar{a}~.$$
The {\it Wu classes} of a chain bundle $(B,\beta)$ are the
$R$-module morphisms
\begin{equation}\label{3c}
v_r(\beta)~:~H_r(B) \to \widehat H^r(\Z_2;R)~;~
x\mapsto \langle\beta_{-2r},x\otimes x\rangle~(r\in \Z)~.
\end{equation}
\indent The {\it universal chain bundle} $(B^R,\beta^R)$ exists for
each $R$.  It is the chain bundle (unique up to equivalence)
characterized by the requirement that the map (\ref{3c}) is an isomorphism
for each $r$.  This implies the more general property that for each
f.g.  free chain complex $C$ the map
\begin{equation} \label{3d}
k_C~:~H_n(C\otimes_R B^R) \to\widehat Q^n(C)~;~
f\mapsto S^{-n}f^\%(\beta)
\end{equation}
is an isomorphism. A cycle $f\in (C\otimes_R B^R)_n$ is a
chain map $f:(B^R)^{-*} \to S^{-n}C$, inducing a morphism
$$S^{-n}f^\%~:~\widehat Q^0((B^R)^{-*}) \to
\widehat Q^0(S^{-n}C)~=~ \widehat{Q}^n(C)~.$$
See Weiss \cite{WeI,WeII} and Ranicki \cite{RAP}.

\subsection{The chain bundle exact sequence and the theorem of Weiss}

For each chain bundle $(B,\beta)$, the map
$J_\beta $ above fits into an exact sequence:
\begin{equation}\label{3e}
\dots\to \widehat Q^{n+1}(B)\overset{H}{\to}
Q_n(B,\beta)\overset{N_\beta}{\to}
 Q^n(B)\overset{J_\beta}{\to} \widehat   Q^n(B)\to \dots
\end{equation}
where the group $Q_n(B,\beta)$ of ``{\it twisted quadratic
structures}'' and the maps $N_\beta$ and $H$ are defined as
follows.

$Q_n(B,\beta) $ is defined as the abelian group of
equivalence classes of pairs $(\phi,\theta)$
(called {\it symmetric structures on $(B,\beta)$})
where $\phi \in(W^\%B)_n$, $\theta \in(\widehat W^\% B)_{n+1}$ satisfy
$$d\phi ~=~0~~,~~d\theta ~=~ J_\beta (\phi)~.$$
The addition is defined by
$$(\phi,\theta) +(\phi',\theta') ~=~ (\phi +\phi',\theta
+\theta'+\xi)~~ \hbox{\rm  where }
\xi_s~=~\phi_0\beta_{s-n+1}\phi'_0~.$$
One says that $(\phi,\theta) $ is equivalent to $(\phi',\theta') $
if there exist
$\zeta\in(W^\%B)_{n+1}$, $\eta\in (\widehat W^\%B)_{n+2}$ such that
$$d\zeta ~=~ \phi'-\phi~~,~~
d\eta ~=~\theta'-\theta+J(\zeta)+(\zeta_0,\phi_0,\phi'_0)^\%(S^n\beta)~.$$
Here $(\zeta_0,\phi_0, \phi'_0)^\%:
(\widehat W^\% B^{-*})_n\to (\widehat W^\%B)_{n+1}$
is the chain homotopy from $\phi_0^{\%}$  to $(\phi')^{\%}_{0} $
induced by $\zeta_0$. (See Ranicki \cite{RAP}, section 3).

The map $H$ is defined by: $H(\theta) =[0,\theta]$.

The map $N_\beta$ is defined by: $N_\beta ([\phi,\theta])= [\phi]$.

When $\beta = 0$, then $Q_n(B,0)=Q_n(B)$ and (\ref{3e}) reduces to
(\ref{3a}).

Recall now from (\cite{RES},p.19,p.39, p.137), the cobordism groups
$L_n(R,\epsilon) $ (resp.  $ L^n(R,\epsilon),\widehat L^n(R,\epsilon
)$) of free $n$-dimensional $\epsilon$-quadratic (resp.  symmetric,
resp.  hyperquadratic) Poincar\'e complexes over $R$, where
$\epsilon=\pm 1$.  These are related by a long exact sequence and a
skew-suspension functor:

\begin{equation}\label{3f}
\begin{CD}
\widehat L^{n+1}(R,\epsilon )   @>H>>   L_n(R,\epsilon )   @>>>
L^n(R,\epsilon )   @>>>  \widehat L^n (R,\epsilon )    \\
@VV \widehat{\overline{S}}^{n+1} V   @VV \overline{S}_n V   @VV
\overline{S}^n V @VV \widehat{\overline{S}}^n V \\
\widehat L^{n+3}(R,-\epsilon )   @>>>   L_{n+2}(R,-\epsilon)
  @>>>L^{n+2}(R,-\epsilon )  @>>>  \widehat L^{n+2}(R,-\epsilon ) .
\end{CD}
\end{equation}
$\overline{S}_n$ is an isomorphism for all $n$, and $L_n(R,1)$ is
the Wall surgery obstruction group, $L_n(R).$ But
$\widehat{\overline{S}}^n$ and $\overline{S}^n$ are not
isomorphisms in general.  Instead, the main result of Weiss
\cite{WeI,WeII} (see also Ranicki \cite{RAP}) identifies the limit
of the maps $\widehat{\overline{S}}^n$ in terms of a functorial
isomorphism:
\begin{equation}\label{3f'}
\lim_{k\to \infty} \widehat L^{n+2k}(R,(-1)^k)
\overset{~\cong~}{\to}Q_n(B^R,\beta^R).
\end{equation}
The skew-suspension maps $\widehat{S}^n$, $S^n$ are isomorphisms
for $1$-dimensional $R$.

\subsection{$ \hbox{\rm UNil} $ and $1$-dimensional rings }

Recall from Definition \ref{1-dim} that a ring $R$ is said to be
$1$-dimensional if it is hereditary and noetherian.

\begin{prop}\label {3h}
For any $1$-dimensional ring $R$ with involution, and any $n\geqslant 0$,
there is a short exact sequence:
$$0\to \hbox{\rm UNil}_n(R;R,R)\to Q_{n+1}(B^{R[x]},\beta^{R[x]})
\to Q_{n+1}(B^R,\beta^R)\to 0$$
\end{prop}
\noindent{\it Proof.}  Following Definition \ref{2.14} set
\begin{equation}\label{3i}
NQ_n(R)~=~ \hbox{\rm ker}\{Q_n(B^{R[x]},\beta^{R[x]})
\to Q_n(B^{R},\beta^{R})\}.
\end{equation}
By Propositions \ref{3g} and \ref{2.1}
$$NL^n(R)~=~L\widetilde{\hbox{\rm Nil}}^n(R)~=~0~\hbox{\rm for
all}~n\geqslant 0~.$$
So by (\ref{3f}) we get a square of isomorphisms, for all $n\geqslant 0$:
\begin{equation}\label{3j}
\begin{CD}
N\widehat L^{n+1}(R,\epsilon)@>~\cong~>> NL_n(R,\epsilon)\\
@V\widehat S^n V ~\cong~ V  @V S_n V ~\cong~ V \\
N\widehat L^{n+3}(R,-\epsilon)@>~\cong~>> NL_{n+2}(R,-\epsilon).
\end{CD}
\end{equation}

By Theorem A, (\ref{3f'}), (\ref{3i}), and (\ref{3j}), for all
$n\geqslant 0$, we have:
$$\begin{array}{ll}
\hbox{\rm UNil}_n(R;R, R) &\cong~ NL_n(R,1)~\cong~ N\widehat
L^{n+1}(R,1)\\[1ex]
&\cong~ \lim_k\, N\widehat L^{n+1+2k}(R,
(-1)^k)~\cong~ NQ_{n+1}(R).
\end{array}$$
This proves (\ref{3h}).
\hfill\qed

\subsection{Rules for calculating $Q_n(C,\gamma)$.}\label{3n}

Our goal, in the light of Proposition \ref {3h}, is to compute
$Q_n(B^A, \beta^A)$, especially when $A=\Z$.  But first we explain
three tools for computing $Q_n(C,\gamma) $ for any chain bundle
$(C,\gamma)$ over any ring with involution $A $.

A) Suppose $(C,\gamma) $ is a chain bundle and
$C\otimes_A C$ is $n$-connected. Then:
$$ Q_i(C,\gamma)~=~0~\hbox{\rm for}~i\leqslant n-1 \quad \hbox{\rm and }
Q^{n+1}(C)\overset{J^{n+1}_\gamma}{\to}\widehat
Q^{n+1}(C)\overset{H^{n+1}}{\to} Q_n(C,\gamma)\to 0
$$
is exact. Moreover, for $i\leqslant n,  J^i_\gamma
=J^i:Q^i(C)\to \widehat Q^i(C), $ and $J^i_\gamma $ is an
isomorphism.

Proof of A): Use the spectral sequence:
$$E^2_{p,q} ~=~
H_p(\Z_2; H_q(C\otimes_A C)) \Rightarrow
H_{p+q}((W^{-*})^{\%} C)~=~Q_{p+q}(C)~.$$
This proves $Q_i(C)=0$, for $i\leqslant n$. Next,
$$J^i_\gamma([\phi])~=~J^i([\phi])-\phi^\%_0([S^i\gamma]) $$
for any $[\phi]\in Q^i(C). $  But if $i\leqslant n,\phi_0 $ is
null homotopic because $[\phi_0]=0\in H_i(C\otimes_A C)$.
Consequently, $\phi_0^\% =0$ and $J^i_\gamma=J^i$ for all $i\leqslant n$.
But by the exact sequence \ref{3a}, it follows that $J^i_\gamma $ is an
isomorphism for all $i\leqslant n$, and $H^{n+1} $ is an
epimorphism. This proves A).

B) Suppose $(C,\gamma)$ is a chain bundle for which the
chain complex $C$ splits as:
$$C~=~\sum\limits_{i=-\infty}^\infty C(i)~.$$
Then
$$\gamma ~=~ \sum\limits_{i=-\infty}^\infty \gamma(i)$$
where $\gamma(i)\in \widehat Q^0(C(i)), $ and the inclusions
$C(i)\to C$ induce a long exact sequence:
\begin{multline}\label {3o}
\dots \to \sum\limits_{i=-\infty}^\infty Q_n(C(i),\gamma(i))\to
Q_n(C,\gamma)\\
\to \sum\limits_{i<j}H_n(C(i)\otimes C(j)) \to
\sum\limits_{i=-\infty}^\infty Q_{n-1}(C(i),\gamma(i)) \to \dots ~.
\end{multline}

Proof of B): On general principles
$$\widehat Q^n(C) ~=~\sum\limits_{i}\widehat Q^n(C(i))$$
and
$$ Q^n(\sum\limits_{i}C(i))~=~ \sum\limits_{i}Q^n(C(i))\oplus
\sum\limits_{i<j}H_n(C(i)\otimes C(j)).$$
Therefore, B) is a consequence of a diagram chase applied
to the following map of exact sequences obtained from
(\ref{3e}):
$$
\begin{CD}
  \sum\limits_{i}
Q_n(C(i),\gamma(i)) @>>> \sum\limits_{i}Q^n(C(i))  @>\Sigma
J_{\beta(i)}>> \sum\limits_{i}\widehat Q^{n}(C(i))  @>>> \\
  @VVV   @VVV   @VV~\cong~ V  \\
   Q_n(C,\gamma) @>>>Q^n(C)  @>{J_\beta}>>  \widehat Q^n(C)
@>>>.
\end{CD}
$$

\

C) Suppose the chain complex $ C $ is concentrated in degrees $\leqslant n$.
Then $ Q^k(C)=0 $ if $k>2n$.  If, in addition, $H_n(C)=0$, then
$Q^{2n}(C) =0 $ as well.

The proof of C. is straightforward from the definition of
$W^\%C$.

The $A$-modules $\widehat{H}^r(\Z_2;A)$ ($r=0,1$) will be said to
be {\it $k$-dimensional} if they admit $k$-dimensional f.g. free
$A$-module resolutions. In the next two subsections we compute
$Q_n(B^A,\beta^A)$ for $A$ with $2A=0$ and $k=0,1$.

\subsection{$Q_n(B^A,\beta^A)$ for $0$-dimensional
$\widehat{H}^*(\Z_2;A)$.  }\label{3o'}

Throughout this section we  suppose $2A =0$, the involution on $A$
is trivial (and consequently $A$ is commutative), and that
$\widehat H^r(\Z_2;A)$ is a f.g. free $A$-module for each $r$.

This occurs, for example, when $A=\mathbb{F}$ or $\mathbb{F}[x]$,
where $\mathbb{F}$ is a perfect field of characteristic 2.

The Frobenius map
$$\psi^2~:~A\to A~;~a \mapsto \psi^2(a)=a^2~.$$
is a ring homomorphism which makes the target copy of $A$ a
module over the source copy of $A$. We denote the target copy
$A$-module as $A'$; thus $A'$ is the additive group of $A$ with
$A$ acting by
$$A \times A' \to A'~;~(a,x) \mapsto a^2x$$
and there is defined an $A$-module isomorphism
$$A' \to \widehat H^r(\Z_2;A)~;~x \mapsto x~.$$

In this case one can easily construct the universal chain
bundle $(B^A,\beta^A)$ for $A$ with
$$d~=~0~:~(B^A)_r~=~A' \to (B^A)_{r-1}~=~A'~.$$
The 0-cycle of $\widehat{W}^{\%}B^{-*}$
$$\beta~=~\sum\limits_r\beta^{-2r}\in
(\widehat{W}^{\%}B^{-*})_0~=~\sum\limits_r
({\rm Hom}_{\Z[\Z_2]}(\widehat{W},B^{-*}_r\otimes_RB^{-*}_r))_0$$
is obtained as follows. Here and below we view $B_r$ as a chain
complex concentrated in degree $r$. Its dual chain complex, $B^{-*}_r$,
concentrated in degree $-r$, consists of $B^r={\rm Hom}_A(B_r,A)$.

Let $x_1\dots x_k$ be a basis of $A'$ over $A$.  Let $x^1\dots x^k$ be
the dual basis.  Write $\underline{x}_i$ for the element $x_i$, viewed
as a member of the ring $A$.  Note that $B^r\otimes_AB^r$ is the
$A$-module of bilinear forms on $B_r$ with values in $A$,
which is canonically identified with
$$(\widehat{W}^{\%}B^{-*}_r)_0~=~
{\rm
Hom}_{\Z[\Z_2]}(\widehat{W}_{-2r},B^{-*}_r\otimes_RB^{-*}_r)~.$$
Therefore the elements $x^i \otimes x^i$,
$\underline{x}_i(x^i\otimes x^i)$ and
$$\beta_{-2r}~:=~\sum\limits_{i=1}^k \underline{x}_i(x^i\otimes x^i)$$
are 0-cycles in $\widehat{W}^{\%}B^{-*}_r$, and bilinear forms
on $B_r$. The matrix of the symmetric bilinear form $\beta_{-2r}$ is diagonal:
$$
\begin{bmatrix}
\underline{x}_1&0&0&\dots\\
0&\underline{x}_2&0&\dots\\
0&0&\underline{x}_3&\dots\\
\vdots&\vdots&\vdots&\ddots
\end{bmatrix}
$$
It follows that $\widehat v_r:H_r(B)\to A'$ is the identity
map. So $(B,\beta)$ is universal. Inclusion induces a map
of chain bundles,
$(B_r,\beta_{-2r})\overset{\iota_r}{\to}(B,\beta)$.

\begin{lem}\label{3p}
Assume $2A=0$, the involution on $A$ is trivial, and $A'$ is free
and finitely generated over $A$. With notation as above, the map
$\iota_r:Q_*(B_r,\beta_{-2r})\to Q_*(B,\beta)$, and the exact
sequence (\ref{3e})\  for $(B_r,\beta_{-2r})$, combine to give an
exact sequence for each r:
\begin{equation}
0\to Q_{2r}(B,\beta)\to
Q^{2r}(B_r)\overset{J_{\beta_{-2r}}}{\to}\widehat Q^{2r}(B_r)\to
Q_{2r-1}(B,\beta)\to 0.
\end{equation}
\end{lem}

\noindent{\it Proof.}
By (\ref{3d}) we have an isomorphism $B_r\otimes
B_{n-r}\overset{k_{B_r}}{~\cong~} \widehat Q^n(B_r). $
 By
\ref{3n}. A), we have $Q_n(B_s,\beta_{-2s})~=~0 \;\hbox{\rm for}
\;n<2s-1 $.
Therefore (\ref{3o}) can be written:
\begin{multline}\label{3q}
\sum\limits_{s\leqslant r}\widehat Q^{2r+1}(B_s)\to \sum\limits_{s\leqslant r}
Q_{2r}(B_s,\beta_{-2s})\to Q_{2r}(B,\beta)\to \sum\limits_{s<
r}\widehat Q^{2r-1}(B_s)\to \\
\sum\limits_{s\leqslant r}
Q_{2r-1}(B_s,\beta_{-2s})\to Q_{2r-1}(B,\beta)\to \sum\limits_{s<
r}\widehat Q^{2r}(B_s)\to Q_{2r-2}(B_s,\beta_{-2s})
\end{multline}
Now, for dimensional reasons, if $n>2s,\; Q^n(B_s)=0, $
and so $\widehat Q^{n+1}B_s)\overset{H }\to
Q_n(B_x,\beta_{-2s})$ is an isomorphism. So (\ref {3q})
reduces to two pieces:
\begin{multline}\label{3r}
\widehat Q^{2r+1}(B_r)\overset{H}{\to}
Q_{2r}(B_r,\beta_{-2r})\to Q_{2r}(B,\beta)\to 0\\
Q_{2r-1}(B_r,\beta_{-2r})\overset{\iota_r}{~\cong~}
Q_{2r}(B,\beta).
\end{multline}
Now apply the exact sequence (\ref{3e}) and Rule \ref{3n} A
to $B_r$ to get:
$$
0\to \hbox{\rm coker}(H_{\beta_{-2r}})\longrightarrow
Q^{2r}(B_r)\overset{J_{\beta_{-2r}}}{\to}\widehat Q^{2r}(B_r)\to
Q_{2r-1}(B_r,\beta_{-2r})\to 0,
$$
which, together with (\ref{3r}) implies Lemma \ref{3p}.
\hfill\qed

 We now restrict ourselves to the case when $A=\mathbb{F}[x]$
where
$\mathbb{F}$ is a perfect field of characteristic 2.  Then
$A'$ is free of rank 2 over $A$, generated by $1$ and $x$.
Since
$B_r= A' $ for all $r$, the abelian group $Q^{2r}(B_r)$ can
be identified with the additive group, $\hbox{\rm Sym}_2(A), $ of
$2\times 2$ symmetric matrices over $A$. The $A$-module $\widehat
Q^{2r}(B_r)$ can be identified with $\hbox{\rm Sym}_2(A)/\hbox{\rm Quad}_2(A)$
where  $\hbox{\rm Quad}_2(A)$ denotes the matrices of the form
$M+M^t$. The map
$J_{{\beta}_{2r}}:\hbox{\rm Sym}_2(A)\to \hbox{\rm Sym}_2(A)/\hbox{\rm
Quad}_2(A)$ then has the form:
$$\begin{array}{ll}
J_\beta
\begin{bmatrix}
a&b\\
b&d
\end{bmatrix}
&=~\begin{bmatrix}
a&b\\
b&d
\end{bmatrix}
\begin{bmatrix}
1&0\\
0&x
\end{bmatrix}
\begin{bmatrix}
a&b\\
b&d
\end{bmatrix}
-
\begin{bmatrix}
a&b\\
b&d
\end{bmatrix}\\[5ex]
&=~\begin{bmatrix}
a^2+a+xb^2&*\\
*&b^2+d+xd^2
\end{bmatrix}~.
\end{array}$$
We intend to show that the kernel and cokernel of $J_\beta$
can be identified with the kernel and cokernel of the map
$\psi^2-1:A\to A$.

We have two inclusion maps $A\overset{\iota}{\to}\hbox{\rm Sym}_2(A)$,
and $A\overset{\iota'}{\to}\hbox{\rm Sym}_2(A)/\hbox{\rm Quad}_2(A)$,
both of the form:
$$ a\to
\begin{bmatrix}
a&0\\
0&0
\end{bmatrix}
$$
Denote the images of these two maps as
$X,X'$. Note that $(J_\beta) \iota =\iota' (\psi^2-1).  $

We use the following easily proved lemma:

\begin{lem} \label{3r'}
Suppose $X,X'$ are subgroups of two abelian groups $Y,Y'$.  Suppose
$j:Y\to Y'$ is a homomorphism such that $j(X)\subseteq X'$, and the
induced map $\tilde{j}:Y/X\to Y'/X'$ is an isomorphism.  Set
$k=j|X:X\to X'$.  Then $\hbox{\rm ker}(k) = \hbox{\rm ker}(j)$, and the
inclusion $X'\to Y'$ induces an isomorphism
$$\iota: \hbox{\rm coker}(k)~\cong~\hbox{\rm coker}(j)~.$$
\end{lem}

We want to apply this lemma when $X,X'$ are as mentioned
earlier and the role of $j: Y \to Y'$  is played by
$$J_{\beta}~:~\hbox{\rm Sym}_2(A) \to \hbox{\rm Sym}_2(A)/\hbox{\rm
Quad}_2(A)~.$$
This means we must first check that
$\tilde{j}$ is an isomorphism.  In other words,  we must check
that each element $p\in \mathbb{F}[x]$ can be written in
one and only one way in the form $b^2+d+xd^2$ where $b,d\in
\mathbb{F}[x]$.

 Write
$$p ~=~\sum\limits_{j=0}^{2n+1}a_jx^j,\quad b~=~ \sum\limits_i b_ix^i,
\quad d~=~\sum\limits_i d_ix^i~.$$
Then:
$$b^2+xd^2+d ~=~ \sum\limits_i(b_i^2+d_{2i})x^{2i}
+\sum\limits_i(d_i^2+d_{2i+1})x^{2i+1}.$$
Therefore the equation $p=b^2+xd^2+d$ reduces to equations,
$$ d_i^2+d_{2i+1}~=~a_{2i+1};\quad b_i^2+d_{2i}~=~a_{2i}.$$
One solves these recursively  for $ d_i $ and $b_i $,
working from higher to lower indices. Note that the first
equation implies that
$d_i=0$ for all $i>n$. Therefore recursively,
the equations
$$d_i^2~=~d_{2i+1}+a_{2i+1}$$
specify $d$.  Then the equations
$$b_i^2 ~=~ d_{2i}+a_{2i}$$
specify $b$.  Here we use that $\mathbb{F} $ is perfect.  Therefore
$\tilde{j}$ is an isomorphism.

Applying the lemma, we conclude that if $A= \mathbb{F}[x]$ then
\begin{equation}\label{3s}
\hbox{\rm ker}(\psi^2-1)\overset{\iota}{~\cong~}  \hbox{\rm ker} J_\beta
\;;\quad  \hbox{\rm coker}(\psi^2-1)\overset{\tilde{\iota}}{~\cong~}
\hbox{\rm coker}(J_\beta)~.
\end{equation}
 The map $\tilde{\iota}$ is induced by $A\overset{\iota'}{\to}
\hbox{\rm Sym}_2(A)/\hbox{\rm Quad}_2(A)$.

Note that if $A=\mathbb{F}_2[x]$, then
$\hbox{\rm ker} (\psi^2-1)=\mathbb{F}_2 $ and the cokernel of
$A\overset{\psi^2-1}{\to}A$ can
be identified with the vector space
$\{\sum\limits_i a_ix^i\,\vert\, a_{2i}=0 \;\hbox{\rm for}\;i>0\}$.

Summarizing, we have a confirmation of the
calculation of Connolly and Ko\'zniewski \cite{CK}:
\begin{thm}: For all $k$, we have  :
$$\begin{array}{l}
\hbox{\rm UNil}_{2k+1}(\mathbb{F}_2;\mathbb{F}_2,\mathbb{F}_2) ~=~0~,\\
\hbox{\rm UNil}_{2k}(\mathbb{F}_2;\mathbb{F}_2,\mathbb{F}_2)~\cong~
{\rm  coker}(\mathbb{F}_2[x]/\mathbb{F}_2\overset{\psi^2-1}{
\longrightarrow}
\mathbb{F}_2[x]/\mathbb{F}_2)\\[1ex]
\hphantom{\hbox{\rm UNil}_{2k}(\mathbb{F}_2;\mathbb{F}_2,\mathbb{F}_2)~}
\cong~  \{\sum\limits_i a_ix^i:a_{2i}=0 \;\hbox{\rm for}\;i\geqslant 0,
a_i \in \mathbb{F}_{2}\}
\end{array}$$
\end{thm}
\noindent{\it Proof.} This is a consequence of Corollary \ref{3h},
Lemma \ref{3p} and  (\ref{3s}).
\hfill\qed

\subsection {$Q_n(B^A,\beta^A)$ for $1$-dimensional $\widehat{H}^*(\Z_2;A)$. }

In this subsection we deal with a ring $A$ whose universal chain
bundle $(B^A,\beta^A)$ satisfies:
\begin{equation} \label{3t}
\hbox{\rm For all $i$,}~
B^A_{2i}\overset{d}{\to} B^A_{2i-1}~
\hbox{\rm is zero;}~ B^A_{2i+1}\overset{d}{\to}B^A_{2i}~
\hbox{\rm is injective}
\end{equation}
with $B^A_r$ f.g. free $A$-modules. Thus $\widehat{H}^0(\Z_2;A)$
has a 1-dimensional f.g. free $A$-module resolution
$$0 \to B^A_{2i+1}\overset{d}{\to}B^A_{2i} \to \widehat{H}^0(\Z_2;A) \to 0$$
and $\widehat{H}^1(\Z_2;A)=0$.
(We shall see that this holds for $A=\Z$ or $\Z[x]$.
The point is that  Corollary \ref{3h}
reduces the calculation of
$\hbox{\rm UNil}_*(\mathbb{Z};\mathbb{Z},\mathbb{Z}) $ to that of
$Q_*(B^A,\beta^A) $ for such rings $A$).

We clearly have:
$$
(B^A,\beta^A)~=~ \sum\limits_{i=-\infty}^{\infty}(B^A(i),
\beta^A(i)),
\hbox{\rm  where $B^A(i) $ is: }\dots 0\to
B^A_{2i+1}\overset{d}{\to} B^A_{2i} \to 0\dots .
$$
We first relate $ Q_n(B^A,\beta^A) $  to
$Q_n(B^A(0),\beta^A(0))$,  for $n=-1,0,1,2,$   by analyzing
the exact sequence   (\ref{3o}), of
the above direct sum splitting.
By (\ref{3n} A), we have:
$$  \sum\limits_{i=-\infty}^{ \infty}Q_m(B^A(i),\beta^A(i)) ~=~
\sum\limits_{i\leqslant \frac{m+1}{4}}Q_m(B^A(i),\beta^A(i)).
$$
Next, because of (\ref{3d}), and dimensional reasons, we
have
$$\begin{array}{ll}
\sum\limits_{i<j}H_m(B^A(i)\otimes B^A(j)) &=~
\sum\limits_{2i<[\frac{m}{2}]}H_m(B^A(i)\otimes
B^A([\dfrac{m}{2}]-i))\\[1ex]
&=~\sum\limits_{2i<[\frac{m}{2}]}H_m(B^A(i)\otimes B^A)\\[1ex]
&=~\sum\limits_{i<\frac{1}{2}[\frac{m}{2}]}\widehat{Q}^m(B^A(i))~.
\end{array}$$
But, by (\ref{3e}) and (\ref{3n}) C),  the map
$\widehat{Q}^{m+1}(B^A(i))\to Q_m(B^A(i),\beta^A(i))$ is an isomorphism if
$i\leqslant \dfrac{m-2}{4} $.
Therefore, after we remove isomorphic direct summands from the
exact sequence (\ref{3o}), it reduces to the much simpler
long exact sequence:
$$\begin{array}{l}
\dots \to\sum\limits_{\frac{m-2}{4}
< i \leqslant \frac{m+1}{4}}Q_m(B^A(i),\beta^A(i))
\to Q_m(B^A,\beta^A)\\
\hskip100pt \to \sum\limits_{\frac{m-3}{4} <i<
\frac{1}{2}[\frac{m}{2}]}\widehat Q^{m}(B^A(i))\to\dots.
\end{array}$$
So,  we get:
\begin{align}\label{3u}
Q_m(B^A,\beta^A) &\overset{~\cong~}{\to}Q_m(B^A(0),
\beta^A(0))
\hbox{\rm \qquad for}~ m=-1,0,~\hbox{\rm and: }\\
Q_1(B^A,\beta^A) &~=~
\hbox{\rm ker}\{Q^1(B^A(0))\overset{J^1_{\beta^A(0)}}{\longrightarrow}\widehat
Q^1(B^A(0))\}\notag \\
Q_2(B^A,\beta^A) &~=~
\hbox{\rm im}\{Q^2(B^A(0))\overset{J^2_{\beta^A(0)}}{\longrightarrow}\widehat
Q^2(B^A(0))\}~=~0~ \hbox{\rm  by}~~(\ref{3n})\, C),\notag
\end{align}
whenever $(B^A,\beta^A)$ is the universal chain bundle of
$A$, and $(B^A,\beta^A)$ satisfies (\ref{3t}).

\

Next we show that (\ref{3t}) holds when $A=\Z$ or $\Z[x]$.

\subsubsection { The construction of $(B^A,\beta^A)$
for certain rings $A$. }\label{3u'}

Suppose $A$ is a commutative ring with no elements of order
2, and trivial involution. Write
$$A_2~=~A/2A~.$$
Therefore $\widehat H^1(\Z_2; A)=0$,  and $\widehat H^0(\Z_2;A) =A_2'$,
by which we mean the abelian group $A_2$,
equipped with the $A$-module structure:
$$ A\times A_2 \to A_2;\; (a, x)\mapsto (a^2x).
$$
Suppose further that there are elements $ x_1, x_2,\dots ,
x_r \in A, r>0,$ such that,
$$
0\to A^r\overset{\times 2}{\to} A^r
\overset{j}{\to }A_2'\to 0
$$
is exact, where
$$j~:~A^r \to A_2'~;~(a_1, a_2,\dots, a_r) \mapsto  a_1^2x_1+
a_2^2x_2 +\dots +a_r^2x_r~.$$
(For example, if $A=\Z$ then $r=1, x_1 = 1$, while if
$A=\Z[x]$ then $r=2, x_1=1, x_2=x$).

We show here how to construct the universal chain bundle
$(B,\beta)$ for $A$, so that (\ref{3t}) holds.

First we construct $B$. For all $i$, we define:
\begin{align}\label{3v}
B_i &~=~ A^r\\
B_{2i}&\overset{d=0}{\longrightarrow}B_{2i-1}\notag\\
B_{2i+1}~=~A^r&\overset{d~=~\times
2}{\longrightarrow}A^r~=~B_{2i}.\notag
\end{align}
Next let $X\in M_r(A)$ be the diagonal matrix,
\begin{equation*}
X~=~
\begin{pmatrix}
x_1&0&\dots &0\\

0&x_2&\dots &0\\
\vdots&\vdots&\ddots&\vdots\\
0&0&\dots&x_r
\end{pmatrix}.
\end{equation*}
We define $\beta=\{\beta_{-i}\in(B^{-*}\otimes B^{-*})_i\}$ by:
\begin{align}\label{3w} \beta_{-4i}&~=~X\in M_r(A) ~=~
(B_{2i}\otimes B_{2i})^*\notag\\
\beta_{-4i-1}&~=~(\delta\otimes 1)\beta_{-4i}\notag\\
\beta_{-4i-2}&~=~-\dfrac{1}{2}(\delta\otimes \delta)\beta_{-4i}~
\hbox{\rm  for all } i.
\end{align}

Here $\delta:B^{-*}_0\to B^{-*}_{-1}$ is the coboundary
homomorphism.

 As in (\ref{3o'}), the map $\widehat
v_{2i}:H_{2i}(B)\to A_2'$ is an isomorphism for all $i$, and
so $(B,\beta)$ is the universal chain  bundle for $A$.

We can now apply the calculation (\ref{3u}) to the
computation of $Q_n(B^A,\beta^A)$, when $A=\Z$ or
$\Z[x]$. Specifically, (\ref{3u}) and Corollary  \ref{3h}
give us the split short exact sequence if $n=0$ or $-1$:
\begin{equation}\label{3w'}
0\to \hbox{\rm UNil}_{n-1}(\Z;\Z,\Z)\to Q_n (B^{\Z[x]}(0),
\beta^{\Z[x]}(0))\overset{\eta_*}{\to}
Q_n(B^\Z(0),\beta^\Z (0))\to 0 ,
\end{equation}
where $\eta:\Z[x]\to\Z$ is the augmentation map.

To simplify things further we define three families of
groups, $K_n, C_n, I_n $, by the exactness of the following
three split sequences:
$$\begin{array}{l}
0\to K_n\to \hbox{\rm ker}
(J_{\beta(0)}^n(\Z[x]))\overset{\eta_*}{\to}
\hbox{\rm ker} (J_{\beta(0)}^n(\Z))\to 0\\[1ex]
0\to C_n \to \hbox{\rm coker} (J_{\beta(0)}^{n+1}(\Z[x]))
\overset{\eta_*}{\to}\hbox{\rm coker} (J_{\beta(0)}^{n+1}(\Z))\to
0\\[1ex]
0\to I_{n+1} \to \hbox{\rm im} (J_{\beta(0)}^{n+1}(\Z[x]))
\overset{\eta_*}{\to}\hbox{\rm im} (J_{\beta(0)}^{n+1}(\Z))\to 0
\end{array}$$
We next claim there is an isomorphism:
\begin{equation}\label{3x_0}
\hbox{\rm UNil}_{-2}(\Z;\Z,\Z) ~\cong~ C_{-1}.
\end{equation}
To see this, note that
$ Q_n(B^A(0))=0$ for dimensional reasons if $ n \leqslant -1$.
Also, by (\ref{3n}) \,A),
$$
J^{-1}_{\beta^A(0)} ~=~ J^{-1}:Q^{-1}(B^A(0))\to \widehat
Q^{-1}(B^A(0)),
$$
 which is a monomorphism by (\ref{3a}). This
implies that
$$ Q_{-1}(B^A(0),\beta^A(0))~\cong~
\hbox{\rm coker}(J^0_{\beta^A(0)})~.$$
Therefore (\ref{3w'}) simplifies when $n=-1$,  to (\ref{3x_0}).

Now  \ref{3u}, Corollary  \ref{3h}, and the exact
sequence (\ref{3e}) for $(B^A(0),\beta^A(0))$ (when $A=\Z,
\Z[x]$)  yield the following calculations:
\begin{align}\label{3x}
&\hbox{\rm UNil}_0(\Z;\Z,\Z) ~\cong~ K_1\\
&\hbox{\rm UNil}_1(\Z;\Z,\Z) ~\cong~ I_2 \notag
\\ 0\to C_0 \to &\hbox{\rm UNil}_{-1}(\Z;\Z,
\Z)\to  K_0 \to 0\notag\\
C_{-1}~\cong~ &\hbox{\rm UNil} _{-2}(\Z;\Z,
\Z).\notag
\end{align}
Therefore our goal is to calculate $C_0, C_{-1}, K_0, $
and $K_1$. This is done in the next two subsections.

\subsubsection{ Calculation of $Q^n(B^A(0))$ and
$\widehat Q^n(B^A(0))$.}

Recall from Ranicki \cite{RATI,RATII} that for any ring with
involution $A$ and for any $A$-module chain complex $C$ an element
$\phi\in(\widehat W^\%C)_n$ is specified by the sequence of
elements $(\dots ,\phi_{-1},\phi_0,\phi_1,\dots)$ of $C \otimes_A
C$ defined by
$$\phi_i~=~\phi(e_i)\in(C\otimes_A C)_{n+i}~(i \in \Z)$$
where $e_i\in \widehat W_i$  is the standard basis element. Likewise,
an element $\phi\in (W^\%C)_n$ is specified by a sequence
$(\phi_0,\phi_1,\dots)$, with  $\phi_i=\phi(e_i)$.

For the rest of this section we assume $A$ is a ring satisfying the
hypotheses at the beginning of section \ref{3u'}.

Let $t:M_r(A)\to M_r(A)$ be the transpose map and define
$$\begin{array}{l}
\hbox{\rm Sym}_r(A)~=~\hbox{\rm ker}(1-t:M_r(A)\to M_r(A))\\[1ex]
\hbox{\rm Quad}_r(A) ~=~ \hbox{\rm im}(1+t:M_r(A)\to M_r(A))
\end{array}$$
Note that $B^A(0)$ is the algebraic mapping cone ${\mathcal C}(f)$
of the map $f:C\to D$, where $C=D=A^r$ is concentrated in degree
0, and $f=\times 2: A^r\to A^r$.  Therefore, for all $m$:
$$
\widehat Q^{2m}(C)~=~\widehat Q^{2m}(D)~\cong~
\hbox{\rm Sym}_r(A)/\hbox{\rm Quad}_r(A):  [\phi]\mapsto [\phi_{-2m}]
$$
because $\phi_{-2m}\in A^r\otimes A^r = M_r(A) $ must be in
the kernel of $1-T$, for all $2m$-cycles $\phi \in
(\widehat W^\%D)_{2m}$.

 Also,
$Q^{2m+1}(C)=\widehat Q^{2m+1}(D) =0$ for all $m$. Since the
induced map, $f^\%:\widehat Q^m(C)\to
\widehat Q^m(D)$ is multiplication by 4, we see $f^\%=0$. So the
sequence:
$$
0\to \widehat Q^m(D)\to \widehat Q^m(B^A(0))\to \widehat Q^m(\Sigma C)\to 0
$$
is exact for all $m$.

If $m=1$ the composite isomorphism,
$$\widehat Q^1(B^A(0))\overset{~\cong~}{\to}\widehat Q^1(\Sigma
C)\overset{~\cong~}{\to}\dfrac{\hbox{\rm Sym}_r(A)}{\hbox{\rm Quad}_r(A)} $$
is written as
\begin{equation}\label{3x_1}
\widehat
Q^1(B^A(0))\overset{\beta^1}{~\cong~}
\dfrac{\hbox{\rm Sym}_r(A)}{\hbox{\rm Quad}_r(A)}: \quad \beta^1([(
\phi_{-1},\phi_{0},\phi_1)]) ~=~ [\phi_1]
\end{equation}

If $m=0$ we write the inverse of the composite isomorphism
$$\dfrac{\hbox{\rm Sym}_r(A)}{\hbox{\rm Quad}_r(A)}~\cong~ \widehat
Q^0(D)~\cong~ \widehat Q^0(B^A(0))$$
as:
\begin{equation}\label{3x_2}
\widehat
Q^0(B^A(0))\overset{\beta^0}{~\cong~}
\dfrac{\hbox{\rm Sym}_r(A)}{\hbox{\rm Quad}_r(A)}: \quad \beta^0([(
\phi_{0},\phi_{1},\phi_2)]) ~=~ [\phi_0]~.
\end{equation}

The calculation of $Q^m(B^A(0))$  requires more work.

Following Ranicki \cite{RATI,RATII} we define $Q^m(f)$ as the
$m$-th homology group of the mapping cone of $f^\%$:
$$Q^m(f)~=~ H_m(f^\%:W^\%C\to W^\%D)~,$$
for any chain map $f:C \to D$ of free $A$-module
chain complexes. We also write $\mathcal {C}(f)$ for the
mapping cone of such $f$, and we write $g:D\to \mathcal {C}(f)
$ for the inclusion. The symmetrization map
\begin{equation*}
H_m(C\otimes_A C)\to Q^m(C)~;~\theta\mapsto\{\phi_s=\begin{cases}
(1+T)\theta &\hbox{\rm if $s=0$}\\
0&\hbox{\rm if $s\geqslant 1$}
\end{cases}
\}
\end{equation*}
fits into a natural transformation of exact sequences:
$$
\begin{CD}
H_m(C\otimes_A C) @>f>>H_m(D\otimes
C)@>g>>H_m(\mathcal {C}(f)\otimes C)@>>>H_{m-1}(C\otimes_A C)\\
@V(1+T)VV  @V(1+T)f VV   @V(1+T)f VV  @V(1+T)VV \\
Q^m(C) @>f^\%>> Q^m(D)@>>>Q^m(f)@>>>Q^{m-1}(C)
\end{CD}
$$
This leads to a further exact sequence relating  $Q^m(f)$
to $Q^m(\mathcal {C}(f))$:
$$
\dots \to Q^{m+1}(\mathcal {C}(f))\to H_m(\mathcal {C}(f)\otimes
C)\overset{(1+T)f}{\longrightarrow}Q^m(f)\to
Q^m(\mathcal {C}(f))\to \dots
$$
Now in the case at hand (where $C=D=A^r$, and $\mathcal {C}(f)=B^A(0)$),
we have
$$Q^{m}(C)~=~Q^{m}(D)~=~
\begin{cases}
\hbox{\rm Sym}_r(A),~\hbox{\rm  if}~m=0\\[1ex]
\dfrac{\hbox{\rm Sym}_r(A)}{\hbox{\rm Quad}_r(A)},~
\hbox{\rm if $m$ is even and $m<0$} \\[1ex]
0,~\hbox{\rm in all other cases }.
\end{cases}
$$
But $f^\%$ is multiplication by 4. Thus
$$\begin{array}{l}
Q^{0}(f)~=~\dfrac{\displaystyle\hbox{\rm
Sym}_r(A)}{\displaystyle4\hbox{\rm Sym}_r(A)}~,\\[2ex]
Q^{2m}(f)~=~\dfrac{\displaystyle\hbox{\rm
Sym}_r(A)}{\displaystyle\hbox{\rm Quad}_r(A)}~
(m< 0)~,\\[2ex]
Q^k(f)~=~0~\hbox{\rm for all other}~k~.
\end{array}$$
So from the above exact sequence, we extract the following diagram with
exact rows:
$$
\begin{CD}
H_0(\mathcal {C}(f)\otimes C) @>{(1+T)f}>>
Q^0(f)@>>>Q^0(B^A(0))@>>>0\\
@VV~\cong~ V  @VV~\cong~ V  @VV ~=~ V\\
\displaystyle{\dfrac{M_r(A)}{2 M_r(A)}}
@>2(1+t)>>\displaystyle{\dfrac{\hbox{\rm Sym}_r(A)}{4\hbox{\rm Sym}_r(A)}}
@>\alpha>> Q^0(B^A(0)) @>>>0
\end{CD}
$$
Therefore $\alpha $  induces an isomorphism:
\begin{equation}\label{3y}
\dfrac{\hbox{\rm Sym}_r(A)}{2 \hbox{\rm Quad}_r(A)} \overset{\alpha^0}{~\cong~}
 Q^0(B^A(0));\quad \alpha^0([M])~=~[(M,0,0)],
\end{equation}
where $(M,0 , 0)$ is a $0$-cycle in $W^\%B^A(0), $ for any
$$M\in \hbox{\rm Sym}_r(A)\subseteq M_r(A) ~=~ A^r\otimes A^r ~=~
(B^A(0)\otimes B^A(0))_0~.$$
\indent
Now $Q^m(B^A(0))=0 $ if $ m\geqslant 2$ by \ref{3n} C). Also by
(\ref{3a}), if $m\leqslant -1$, the map $Q^m(B^A(0))\overset
{J^m}{\to} \widehat Q^m(B^A(0))$ is an isomorphism.

Therefore, we are only left with the calculation of $Q^1(B^A(0))$.
Instead of the above method (which would yield the result) we calculate
this by hand both for its therapeutic value and for its greater
explicitness.  The bottom line will be (\ref{cc}).

For each $M\in M_r(A)$, define
$$\phi^M~=~(\phi^M_0,\phi^M_1)\in( W^\%B^A(0))_1 $$
by:
\begin{align*}
\phi^M_1 &~=~ M \in M_r(A)~=~ A^r\otimes A^r ~=~ B_1\otimes
B_1\notag\\
\phi_0^M &~=~ M\oplus (-M)\in (B_1\otimes B_0)\oplus
(B_0\otimes B_1)
\end{align*}
where $B_i=B^A(0)_i$.

\begin{lem}
If $M\in \hbox{\rm Sym}_r(A)$, then $\phi^M$ is
a 1-cycle in $W^\%B^A(0)$, and the rule $M\mapsto \phi^M$
induces an isomorphism:
\begin{equation}\label{cc}
\alpha^1~:~\dfrac{\hbox{\rm Sym}_r(A)}{2\hbox{\rm Sym}_r(A)}~\cong~Q^1(B^A(0)).
\end{equation}
\end{lem}

\noindent{\it Proof.}
For any $\phi=(\phi_0,\phi_1)\in (W^\%B^A(0))_1$, $\phi=(\phi_0,\phi_1)$
where $\phi_i \in (B^A(0)\otimes_AB^A(0))_{i+1}$. We can write
$$\phi_0~=~ \kappa_1\oplus \kappa_2~,$$
where $\kappa_1\in M_r(A)=B_1\otimes_A B_0, $ and $ \kappa_2\in
M_r(A)=B_0\otimes_AB_1$.  $\phi $ is a 1-cycle if and only if:
$$
1)\;  \p \phi_0 ~=~0;\qquad 2)\; (T-1)\phi_0 ~=~ -\p \phi_1;
\qquad 3)\;(T+1)\phi_1 ~=~0~,$$
where $T:B^A(0)\otimes B^A(0)\to B^A(0)\otimes B^A(0)$ is
the twist chain map
$$T(x\otimes y) ~=~(-1)^{|x||y|}y\otimes x~.$$
These three conditions are equivalent to:
$$\kappa_2 ~=~ -\kappa_1;\qquad (1+t)\kappa_1 ~=~ 2\phi_1;\qquad
t\phi_1~=~ \phi_1~ \hbox{\rm  in } A^r\otimes_A A^r ~=~ M_r(A)~.$$
Here $t$ denotes the transpose map in $M_r(A)$.
Also a cycle $\phi$ as above is a boundary in $W^\%B^A(0)$ if
and only if there is an element $\psi \in B_1\otimes B_1, $
such that $ \kappa_1=2\psi $ in $A^r\otimes A^r=M_r(A)$.
Therefore the map
\begin{equation}\label{5}
Q^1(B^A(0))\to \hbox{\rm Sym}_r(A)/2\hbox{\rm Sym}_r(A)~:~
[\phi]\mapsto \kappa_1 \bmod (2A)
\end{equation}
is an isomorphism.

The above discussion shows that if $M\in \hbox{\rm Sym}_r(A)$, then
$\phi^M$ is a 1-cycle, and if $M\in 2\hbox{\rm Sym}_r(A)$, then
$\phi^M$ is a boundary. Since the map (\ref{5})
obviously sends $\phi^M$ to $M$, the proof is complete.
\hfill\qed

We summarize the calculations of this subsection as follows:
\begin{align}\label{3y'}
\widehat Q^m(B^A(0))&\;~\cong~ \dfrac{\hbox{\rm Sym}_r(A)}{\hbox{\rm
Quad}_r(A)}~
\hbox{\rm   for all } m,\notag \\
 Q^0(B^A(0))&~\cong~
\dfrac{\hbox{\rm Sym}_r(A)}{2\hbox{\rm Quad}_r(A)}\\
Q^1(B^A(0))&~\cong~ \dfrac{\hbox{\rm Sym}_r(A)}{2\hbox{\rm Sym}_r(A)} \notag \\
Q^n(B^A(0))&~=~0~ \hbox{\rm  for } n\geqslant 2 \notag\\
Q^n(B^A(0))&\overset{J^n}{ ~\cong~ }\widehat Q^n(B^A(0))~
\hbox{\rm  if } n \leqslant -1.  \notag
\end{align}

\subsubsection{ The maps $J^0_{\beta(0)}(A),
J^1_{\beta(0)}(A)$ and the groups $C_{-1},C_{0}$  and $K_1, K_0$.}

We first analyze the map $J^0_{\beta(0)}(A):Q^0(B^A(0))\to
\widehat Q^0(B^A(0))$, when $A=\Z$ or $\Z[x]$ using
the isomorphisms of
(\ref{3x_1}),(\ref{3x_2}),(\ref{3y}),(\ref{cc}). By \ref {3b},
$\beta^0\circ J^0_{\beta(0)}(A)\circ \alpha^0$
 sends a matrix $M\in\dfrac{\hbox{\rm Sym}_r(A)}{2\text{Quad}_r(A)}$ to:
$$
\beta^0 (J^0([(M, 0, 0)]))-M^tXM  ~=~
M-MXM\in
\dfrac{\hbox{\rm Sym}_r(A)}{\hbox{\rm Quad}_r(A)}.
$$
In the case when $A=\Z$, so that $r=1$, and $ X=1, $ we have
$\beta^0
J^0_{\beta(0)}(\Z)\alpha^0  $, sending
$a\in \Z_4$ to
$a-a^2 \in \Z_4/2\Z_4=\Z_2$. So $J^0_{\beta(0)}(\Z) = 0. $
Therefore:
$$
\hbox{\rm ker}J^0_{\beta(0)}(\Z)~=~ Q^0(B^{\Z}(0))~\cong~ \Z_4;\quad
\hbox{\rm coker}J^0_{\beta(0)}(\Z)~=~\widehat Q^0(B^\Z(0))~\cong~ \Z_2.
$$
Now we let $A=\Z[x]$.
Set
$$\mathcal {J}^0~=~\beta^0\circ
J^0_{\beta(0)}(\Z[x])\circ\alpha^0:
\hbox{\rm Sym}_2(\Z[x])/2\hbox{\rm Quad}_2(\Z[x])\to \hbox{\rm
Sym}_2(\Z[x])/\hbox{\rm Quad}_2(\Z[x])~.$$

For any
$$\begin{pmatrix}
a&b\\
b&d
\end{pmatrix}
\in \hbox{\rm Sym}_2(\Z[x])/2\hbox{\rm Quad}_2(\Z[x])$$
we compute from the above formula:
\begin{equation}\label{3z}
\mathcal {J}^0
\begin{pmatrix}
a&b\\
b&d
\end{pmatrix}
~=~
\begin{pmatrix}
a-a^2-b^2x&b-ab-bdx\\
b-ab-bdx&d-b^2-d^2x
\end{pmatrix}
\in \dfrac{\hbox{\rm Sym}_2(A)}{\hbox{\rm Quad}_2(A)}
\end{equation}

We want to apply Lemma \ref{3r'} again. Let
$j=\mathcal {J}^0$, and:
$$Y ~=~\dfrac{\hbox{\rm Sym}_r(\Z[x])}{2\hbox{\rm Quad}_r(\Z[x])},\;
Y'~=~\dfrac{\hbox{\rm Sym}_r(\Z[x])}{\hbox{\rm Quad}_r(\Z[x])},\;
X~=~(\Z_4[x])\times (\F_2[x]),\;X'~=~(\F_2[x])
$$

 $X$ and $X'$ include into $Y$ and $Y'$ respectively
by the rules:
$(a,d)\mapsto
\left(
\begin{smallmatrix}
a&0\\
0&2d
\end{smallmatrix}
\right)
$, and
$a \mapsto
\left(
\begin{smallmatrix}
a&0\\
0&0
\end{smallmatrix}
\right)
$.
We first will have to show
 that $Y/X\to Y'/X'$ is an isomorphism.
To this end, we note an isomorphism,
$\F_2[x]\times\F_2[x]~\cong~ Y/i(X), $ defined by:
$(b,d)\mapsto
\left(
\begin{smallmatrix}
0&b\\
b&d
\end{smallmatrix}
\right),
$ and an isomorphism $\F_2[x]~\cong~ Y'/i'(X'), $ given by
: $p\mapsto
\left[
\begin{smallmatrix}
0&0\\
0&p
\end{smallmatrix}
\right]
$.
Therefore the claim that $j$ induces an isomorphism,
$Y/X\to Y'/X'$, amounts to the statement that each
$p\in\F_2[x]$ can be written uniquely in the form,
$p~=~b^2+d+xd^2, $ for some $b, d\in \F_2[x]. $ But this
 was proved already in section \ref{3o'}.

Define
$$k~:~\Z_4[x]\times \F_2[x] \to\F_2[x]~;~
(a,d) \mapsto a-a^2 \bmod 2~.$$
Clearly,
$$\begin{array}{l}
\hbox{\rm ker}(k)~=~\{(a,d)\in \Z_4[x]\times \F_2[x]|
a=a_0+2a_1,\hbox{\rm  for some } a_0\in \Z_4,
a_1\in \Z_4[x]\}~,\\[1ex]
\hbox{\rm coker}(k)~=~\hbox{\rm coker}(\psi^2-1)~.
\end{array}$$
Applying Lemma (\ref{3r'}), we
see that $i$ and $i'$ induce isomorphisms:
$$
\hbox{\rm ker}(k)\overset{\iota}{~\cong~} \hbox{\rm ker}J^0_{\beta(0)}
(\Z[x]);\quad \hbox{\rm coker}(\psi^2-1)\overset{\iota'}{\to}
\hbox{\rm coker}(J^0_{\beta(0)} (\Z[x])).
$$

Also, $\iota(a,d) ~=~ \alpha^0 \left[
\begin{smallmatrix}
a&0\\
0&2d
\end{smallmatrix}
\right]$.

The augmentation map induced by $\eta$
$$Q^0(B^{\Z[x]}(0))\overset{\eta_*}{\to}Q^0(B^\Z(0))$$
sends $\alpha^0
\left[
\begin{smallmatrix}
a&0\\
0&2d
\end{smallmatrix}
\right]
$ to $a_0\in \Z_4, $ the degree zero coefficient of $a$.
The same formula holds as well for
$\eta_*:Q^0(B^{\Z[x]}(0)) \to Q^0(B^\Z(0)).$

Restricting $\eta_*$ to $\hbox{\rm ker}J^0_{\beta(0)} (\Z[x])$, we
get a short exact sequence:
$$
0\to \F_2[x]\times
\F_2[x]\overset{k_2}{\to}\hbox{\rm ker}(J^0_{\beta(0)}(\Z[x]))
\overset{\eta_*}{\to}\hbox{\rm ker}(J^0_{\beta(0)}(\Z))\to 0
$$
where $k_2$ is defined by:
$$k_2(a,d)~=~
\begin{pmatrix}
2a&0\\
0&2d
\end{pmatrix}~.$$
This yields isomorphisms:
\begin{equation}\label{bb}
\F_2[x]\times \F_2[x]\overset{k_2}{~\cong~}K_0,\quad
\hbox{\rm coker}\{(\psi^2-1):\F_2[x]/\F_2\to\F_2[x]/\F_2\}
\overset{k_2'}{~\cong~} C_{-1}.
\end{equation}
Here $ (\psi^2-1):\F_2[x]/\F_2\to\F_2[x]/\F_2  $  \;is the
map induced by $\psi^2-1:\F_2[x]\to\F_2[x], $ and $k_2'$
is induced by $\iota'$.

Now we analyse $J^1_{\beta(0)}(A)$ similarly.
Recall $B^A(0)$ is a chain complex concentrated in degrees
$0$ and $1$: $B_0=A^r;\; B_1= A^r,\;$ and its boundary map
is $\p=\times 2:B_1\to B_0$.

In order to understand the map $J^1_{\beta(0)}(A)$, we
define, for any 1-cycle, $\phi\in (W^\%B^A(0))_1$, another
1-cycle
$$
\gamma^\phi ~=~\phi_0^\%(S^1(\beta(0)))\in (\widehat W^\%B^A(0))_1\;
.
$$
We know $ \gamma^\phi=(\gamma^\phi_{-1},\gamma^\phi_0,
\gamma^\phi_1),$ where
$$\gamma^\phi_i~=~\gamma_i~=~ \tilde
\phi_0\otimes \tilde\phi_0(\beta(0)_{i-1})~.$$
Here $\tilde \phi_0:B^A(0)^{1-*}\to B^A(0)$ is the chain
map whose matrix is $\phi_0\in (B^A(0)\otimes B^A(0))_1. $

We conclude:
$$\begin{array}{l}
\gamma_1~=~ \tilde \phi_0\otimes \tilde\phi_0(X)\in
B_1\otimes B_1,\\[1ex]
\gamma_0~=~(1\otimes \p)\gamma_1\in (B^A(0)\otimes
B^A(0))_1,\\[1ex]
\gamma_{-1}~=~ \dfrac{1}{2}(\p\otimes\p)\gamma_1 \in B_0\otimes B_0.
\end{array}$$

Therefore
$$
J^1_{\beta(0)}(A):Q^1(B^A(0))\to \widehat
Q^1(B^A(0))\hbox{\rm \quad is : } [\phi]\mapsto
J^1([\phi])-[\gamma^\phi].
$$

Set
$$\mathcal {J}^1~=~\beta^1\circ J^1_{\beta(0)}(A)\circ\alpha^1~.$$
We get:
$$\begin{array}{ll}
\mathcal {J}^1(M) &=~ \beta^1(J^1[\phi^M])-[\gamma_1^{\phi^M}]~=~
M-M^tXM\\[1ex]
&=~M-MXM\;\bmod \hbox{\rm Quad}_r(A)
\end{array}$$
for all $M\in \hbox{\rm Sym}_r(A). $ (The formulae for $\mathcal{J}^1$
and $\mathcal {J}^0$ are identical!). Therefore the formula
(\ref{3z}) can also be used for $\mathcal J^1.$
We therefore conclude at once that we have an isomorphism,
induced by $\beta^1$:
\begin{equation}\label{dd}
\hbox{\rm coker}\{(\psi^2-1):\F_2[x]/\F_2\to\F_2[x]/\F_2\}~\cong~
C_{0}
\end{equation}
To compute $K_1$, we note from \ref{3z} that the kernel
of  $J^1_{\beta(0)}(\Z[x])\circ\alpha^1$ is:
$$
\{
\left[
\begin{smallmatrix}
a_0&0\\
0&0
\end{smallmatrix}
\right]
\in \hbox{\rm Sym}_2(\F_2[x]): a_0\in \F_2\}.
$$
Since $\widehat \eta_*
\left[
\begin{smallmatrix}
a_0&0\\
0&0
\end{smallmatrix}
\right]~=~ a_0\in \Z_2,
$
we conclude at once that:
\begin{equation}\label{ee}
K_1 ~=~ 0
\end{equation}

\subsection{ The calculation of $\hbox{\rm UNil}_n(\Z;\Z,\Z)$ for all $n$.}

The results of the last section allow us to prove Theorem B of the
Introduction:

\begin{thm}\label{ff}
There are isomorphisms:
 \begin{align*}
\hbox{\rm UNil}_0(\Z;\Z,\Z)&=~0,\\
\hbox{\rm UNil}_1(\Z;\Z,\Z)&=~0,\\
\hbox{\rm UNil}_2(\Z;\Z,\Z)&\cong~ \hbox{\rm coker}\{ \;(\psi^2-1): \F_2[x]
/\F_2\to
\F_2[x] /\F_2 \;\},
\end{align*}
and an exact sequence:
$$0\to \F_2[x] /\F_2 \overset{(\psi^2-1)}{\longrightarrow} \F_2[x]/\F_2
\to     \hbox{\rm UNil}_{3}(\Z;\Z,\Z)\to \F_2[x]\times \F_2[x] \to 0~.$$
\end{thm}
\noindent{\it Proof.} Note $I_2=0$, by (\ref{3u}).  Therefore (\ref{ee}) and
(\ref{3x}) imply the first two equations at once.  The third equation is
immediate from (\ref{3x}) and (\ref{bb}).  The final exact sequence is
immediate from (\ref{3x}), (\ref{bb}), and (\ref{dd}).  \hfill\qed

See Banagl and Ranicki \cite{BR} and Connolly and Davis \cite{CD}
for further computations.

\end{document}